\documentclass[11pt, reqno]{amsart}
\usepackage[mac]{inputenc}
\theoremstyle{definition}
\newtheorem*{ack}{Acknowledgements}

\newtheorem{conjecture}{Conjecture}[section]
\newtheorem{example}{Example}[section]
\newtheorem{remark}{Remark}[section]
\newtheorem{definition}{Definition}[section]

\theoremstyle{theorem}
\newtheorem{theorem}{Theorem}[section]
\newtheorem{corollary}{Corollary}[section]
\newtheorem{proposition}{Proposition}[section]
\newtheorem{lemma}{Lemma}[section]

\newtheorem{maintheorem}{Theorem}

\usepackage{amssymb, amsmath, mathrsfs}
\usepackage{fourier, xypic}
\usepackage{hyperref}
\usepackage{mathtools}
\usepackage{color}

\usepackage{enumitem}
\makeatletter
\def\namedlabel#1#2{\begingroup
    #2%
    \def\@currentlabel{#2}%
    \phantomsection\label{#1}\endgroup
}
\makeatother

\newcommand{\alg}{\overline{\mathbb{Q}}}

\newcommand{\CC}{\mathbb{C}}
\newcommand{\FF}{\mathbb{F}}
\newcommand{\PP}{\mathbb{P}}
\newcommand{\QQ}{\mathbb{Q}}
\newcommand{\ZZ}{\mathbb{Z}}
\newcommand{\HH}{\mathbb{H}}

\newcommand{\Ec}{\mathcal{E}}
\newcommand{\Oc}{\mathcal{O}}
\newcommand{\Mca}{\mathcal{M}}
\newcommand{\Fca}{\mathcal{F}}
\newcommand{\Lcal}{\mathcal{L}}

\newcommand{\id}{\mathrm{id}}
\newcommand{\Res}{\mathrm{Res}}
\newcommand{\End}{\mathrm{End}}
\newcommand{\Hom}{\mathrm{Hom}}
\newcommand{\per}{\mathrm{per}}
\newcommand{\rk}{\mathrm{rk}}
\newcommand{\CH}{\mathrm{CH}}
\newcommand{\coker}{\mathrm{coker}}
\newcommand{\Gal}{\mathrm{Gal}}
\newcommand{\tr}{\mathrm{tr}}
\newcommand{\cl}{\mathrm{cl}}
\newcommand{\disc}{\mathrm{disc}}
\newcommand{\ch}{\mathrm{ch}}
\newcommand{\Td}{\mathrm{Td}}
\newcommand{\Gr}{\mathrm{Gr}}
\newcommand{\ord}{\mathrm{ord}}
\newcommand{\Tr}{\mathrm{Tr}}

\newcommand{\unit}{\mathbf{1}}

\newcommand{\cyc}{c_{X \mathrm{mod} D}}

\newcommand{\abs}[1]{\lvert#1\rvert}

\title[Periods of Hodge structures and gamma values]{Periods of Hodge structures and \\ special values of the gamma function}

\author{Javier Fres\'an}

\address{ETH Z\"urich, D-MATH, R\"amistrasse 101, CH-8092 Z\"urich, Switzerland}

\email{javier.fresan@math.ethz.ch}

\begin{document}

\maketitle

\begin{abstract} At the end of the 70s, Gross and Deligne conjectured that periods of geometric Hodge structures with multiplication by an abelian number field are products of values of the gamma function at rational arguments, with exponents determined by the Hodge decomposition. We prove an alternating variant of this conjecture for smooth projective varieties acted upon by an automorphism of finite order, thus improving previous results of Maillot and R\"ossler. The proof relies on a product formula for periods of regular singular connections due to Saito and Terasoma.  
\end{abstract}

\tableofcontents

\section{Introduction}

This paper pursues a series of works by Weil \cite{Wei76}, Gross \cite{Gro78}, Shimura \cite{Shi79}, Deligne \cite{Del82}, Anderson \cite{And82}, Colmez \cite{Col93} and, more recently, Maillot and R\"ossler \cite{MR04}, aiming to understand the relations between periods of certain motives over $\alg$ and values of the gamma function at rational arguments. Recall that the latter is defined by the convergent integral 
$$
\Gamma(z)=\int_0^\infty t^{z-1} e^{-t} dt
$$ for $\mathrm{Re}(z)>0$, and extended to a nowhere vanishing meromorphic function on the complex plane, with simple poles at $z=0, -1, -2, \ldots$

\subsection{Euler's formula} \

\noindent
The oldest instance of such a relation is Euler's formula
\begin{equation*}\label{eq-1}
\int_0^1 t^{a-1}(1-t)^{b-1}dt=\frac{\Gamma(a)\Gamma(b)}{\Gamma(a+b)}, \qquad \mathrm{Re}(a), \ \mathrm{Re}(b)>0, 
\end{equation*} the left hand side of which may be interpreted, for rational values of $a$ and $b$, as a period of a Fermat curve. Let us briefly recall how. 

Leaving aside some trivial cases, we may assume $a=r/d$ and $b=s/d$ for integers $d \geq 3$ and $0<r, s<d$ such that $r+s\neq d$. Choose a primitive $d$-th root of unity $\xi$ and denote by $K$ the cyclotomic field $\QQ(\xi)$. Consider the curve $C \subset \PP^2$ defined by the equation $x^d+y^d=z^d$, which we shall regard as a $K$-variety. The group $G=\mu_d(\CC)^3 \slash (\mathrm{diag})$ acts on $C$ by coordinate-wise multiplication, and this action decomposes the algebraic de Rham cohomology into one-dimensional representations; they are spanned by the classes of the differentials $\omega_{r, s}=x^{r-1} y^{s-1} dx$, on which $G$ acts through the characters $\chi_{r, s}(\xi^i, \xi^j)=\xi^{ir+js}.$ Furthermore, the singular homology of the Riemann surface $C_{\CC}^{\textup{an}}$ is a cyclic $\QQ[G]$-module, and there exists a generator\footnote{Explicitly, $\kappa=\gamma-(1, \xi)_\ast \gamma+(\xi, \xi)_\ast\gamma-(\xi, 1)_\ast \gamma$, where $\gamma$ is the path $t \mapsto (t, (1-t^d)^{1/d})$, $t \in [0, 1]$, on the real points of the affine curve $x^d+y^d=1$.} $\kappa$ \cite{Roh78} such that
\begin{equation}\label{eq-2}
\int_\kappa \omega_{r, s}=\frac{1}{d} (1-\xi^r)(1-\xi^s) \int_0^1 t^{a-1} (1-t)^{b-1} dt. 
\end{equation}

From the standpoint of this article, Euler's formula should be understood as an expression for the determinant of periods of local systems on open varieties purely in terms of local monodromies at infinity. To make this more precise, observe that the quotient $\pi: C \to C \slash G$ realizes $C$ as a $G$-cover of $\PP^1$, ramified at $D=\{0, 1, \infty\}$. By transport of structure, $G$ acts on the direct image of the exterior differential $\Oc_C \to \Omega^1_C$, splitting it into a sum of rank one logarithmic connections 
$$\nabla_\chi: \Ec_\chi \to \Ec_\chi \otimes_{\Oc_{\PP^1}} \Omega^1_{\PP^1}(\log D)
$$ indexed by the characters $\chi$ of $G$. The sheaf of analytic horizontal sections of $\Ec_\chi$ restricts to a complex local system $\ker(\nabla^\textup{an}_\chi)$ on the analytic manifold $U_\CC^\textup{an}$ associated to $U=\PP^1 \setminus D$; it is canonically isomorphic to the $\chi$-component of $(\pi_\ast \CC)_{| U_\CC^{\textup{an}}}$. Since this decomposition is already defined over $K$, we get a local system of $K$-vector spaces $V_\chi$ on $U_\CC^\textup{an}$ together with an isomorphism $\rho \colon V_\chi \otimes_K \CC \stackrel{\sim}{\to} \ker(\nabla_\chi^\textup{an}).$ At the level of cohomology, this $K$-structure induces ``period isomorphisms'' 
\begin{equation*}\label{eq-3}
H^j(\rho) \colon H^j(U, DR(\Ec_\chi, \nabla_\chi)) \otimes_K \CC \stackrel{\sim}{\longrightarrow} H^j(U_\CC^\textup{an}, V_\chi) \otimes_K \CC. 
\end{equation*}

If $\chi$ is non-trivial, the cohomology vanishes in degrees $j=0, 2$ and has dimension one for $j=1$, so the alternating determinant of $H^\bullet(\rho)$ with respect to $K$-bases is simply a period; for $\chi=\chi_{r, s}$, it coincides with the integral \eqref{eq-2} up to an element of $K^\times$. Now a straightforward computation shows that the residues of $\nabla_{\chi_{r, s}}$ at the singular points $0,$ $1$ and $\infty$ are $a,$ $b$ and $-a-b$ respectively, the arguments at which the gamma function is evaluated in Euler's formula!

\subsection{Plan} \

\noindent 
Following an idea of Bloch \cite{Blo05}, we shall use a far-reaching generalization of Euler's formula, the Saito-Terasoma theorem \cite{ST97}, to compute the periods of varieties endowed with automorphisms of finite order. This will yield new cases of a conjecture of Gross and Deligne expressing periods of geometric Hodge structures with multiplication by an abelian number field as products of special values of the gamma function. To our knowledge, the present article constitutes the first application of this fundamental result besides the proof in \textit{loc. cit.} that the category of rank one motives associated to algebraic Hecke characters is closed under taking the determinant of cohomology. 

The paper is organized as follows. 

In Sect. \ref{sec-1}, we recall the Gross-Deligne conjecture and state our main results. Sect. \ref{sec-2} is devoted to a self-contained presentation of the Saito-Terasoma theorem, as well as some complements needed for the sequel. In Sect. \ref{sec-3}, we introduce cyclic covers \`a la Esnault-Viehweg and compute their periods using the Saito-Terasoma theorem. From this we derive the main results in Sect. \ref{sec-4}. 

\subsection{Conventions} \

\noindent 
The following notation and terminology are used throughout this paper: 

\begin{itemize}

\item $\alg$ denotes the algebraic closure of $\QQ$ in $\CC$.

\item If $K$ is a subfield of $\alg$, we write $\sim_{K^\times}$ for the equivalence relation on non-zero complex numbers given by $a \sim_{K^\times} b$ if and only if $a/b \in K^\times$. 

\item If $x \in \QQ$, $[x]$ and $\{x\}$ denote the integral and fractional part respectively.

\item By an \textit{algebraic variety} over a field $k$ we mean an integral separated $k$-scheme of finite type. If $X$ is an algebraic variety over a subfield $k$ of $\CC$, we denote by $X_\CC^\textup{an}$ the complex analytic manifold associated to $X \times_k \CC$, by $H^\ast_B(X)$ (resp. $H_{B, c}^\ast(X)$) the singular cohomology $H^\ast(X^\textup{an}_\CC, \QQ)$ (resp. the singular cohomology with compact support) and by $\chi(X)$ the Euler characteristic of $X_\CC^\textup{an}$. 

\item We say that a divisor $D$ on $X$ has \textit{simple normal crossings} if all its irreducible components $(D_i)_{i \in I}$ are smooth and intersect transversally. For each $i \in I$, we denote by $D_i^\circ$ the complement in $D_i$ of the union of all the irreducible components of $D$ distinct from $D_i$. 

\item If $\mathcal{F}$ is a locally free $\Oc_X$-module, $\mathbf{V}(\Fca)$ denotes the associated geometric vector bundle $\mathbf{Spec}(\mathrm{Sym}^\bullet(\Fca^\vee))$.  

\end{itemize}

\section{Main results}\label{sec-1}

In this section, we first introduce periods of geometric Hodge structures with multiplication by an abelian number field and recall the Gross-Deligne conjecture following \cite[\S 4]{Gro78} and \cite[\S 4]{MR04}. After these preliminaries, we state our main results. For a more extensive discussion on special values of the gamma function as periods of motives, as well as their relation to algebraic Hecke characters, we refer the reader to \cite[\S 8]{Del79}, \cite{Sch88} and \cite[Chap. 24]{And04}. 

\subsection{Periods of Hodge structures with abelian multiplication} \

\noindent 
Throughout, $K$ is a subfield of $\alg$ and $F$ a number field embeddable into $K$. 

\begin{definition}
A \textit{Hodge structure with $F$-multiplication} is a pair $(H, \iota)$ consisting of a rational pure Hodge structure $H$ such that $\dim_\QQ H=[F: \QQ],$ and a ring morphism $\iota: F \to \End_{HS}(H)$. 
\end{definition}

To each such pair it is attached a canonical decomposition
\begin{equation*}
H \otimes_\QQ \CC=\bigoplus_{\sigma \in \Hom(F, \CC)} H_\sigma, 
\end{equation*} where $H_\sigma$ denotes the one-dimensional complex vector space $H \otimes_{F, \sigma} \CC$. Since $F$ acts by morphisms of Hodge structures, $H_\sigma$ is contained in a unique summand $H^{p, q}$ of the Hodge decomposition. Let $p: \Hom(F, \CC) \to \ZZ$ be the function assigning $p$ to $\sigma$. If $H$ has weight $m$, then $p(\sigma)+p(\bar{\sigma})=m$. 

In order to define periods, we shall further assume that:
\begin{enumerate}
\item[\namedlabel{h1}{(H1)}]   $H$ is a sub-Hodge structure of the Betti cohomology $H^m_B(X)$ of a smooth projective variety $X$ over $K$.

\item[\namedlabel{h2}{(H2)}] The subspaces $H_\sigma \subset H \otimes_\QQ \CC$ are defined over $K$ for the ``de Rham structure'' steaming from the comparison isomorphism
$$
\rho^m \colon H^m_{dR}(X \slash K) \otimes_K \CC \stackrel{\sim}{\longrightarrow} H^m_B(X) \otimes_\QQ \CC. 
$$
\end{enumerate}

\begin{remark} Conditions \ref{h1} and \ref{h2} automatically hold when $H$ is the Hodge realization of a motive over $K$ with coefficients in $F$, or if $F$ is a CM field and $H$ is polarizable and effective\footnote{In the latter case, $X$ can be chosen to be a product of CM abelian varieties \cite[Thm. 3]{Abd05}.}. 
\end{remark}

Assuming \ref{h1} and \ref{h2}, for each $\sigma \in \Hom(F, \CC)$ there exists $\omega_\sigma \in H^m_{dR}(X \slash K)$ such that $\rho^m(\omega_\sigma \otimes_K 1_\CC)$ is a basis of $H_\sigma$. Let $\gamma \in H^\vee$ be any non-zero element of the dual of the $\QQ$-vector space $H$. Then the non-zero complex number
$$
\langle \omega_\sigma, \gamma \rangle := (\gamma \otimes_\QQ 1_\CC)(\rho^m(\omega_\sigma \otimes_K 1_\CC))
$$ is independent of the choices of $\omega_\sigma$ and $\gamma$ up to multiplication by an element of $K^\times$. Viewing $\gamma$ as a cycle in $H_m(X_\CC^\textup{an}, \QQ)$, it is nothing else than $\int_\gamma \omega_\sigma$. 

\begin{definition}
The period $P_\sigma(H) \in \CC^\times \slash K^\times$ is the class of $\langle \omega_\sigma, \gamma\rangle$ modulo $K^\times$.
\end{definition}

Suppose finally that the number field $F$ is abelian of conductor $d$ and choose embeddings $F \subseteq \QQ(\mu_d) \subset \alg$. Identifying $\Gal(F \slash \QQ)$ with $\Hom(F, \CC)$, one gets a surjective map $\varphi \colon (\ZZ/d)^\times \to \Hom(F, \CC)$. Given $\lambda \in (\ZZ/d)^\times$, we shall write $P_\lambda(H)$ for the corresponding period $P_{\varphi(\lambda)}(H)$. Besides, the composition of $\varphi$ and $p$ gives a function from $(\ZZ/d)^\times$ to $\ZZ$, still denoted by $p$, which now satisfies $p(\lambda)+p(-\lambda)=m$ since $-1 \in (\ZZ/d)^\times$ acts by complex conjugation. By \cite[Lemma 6.12]{Del77}, there exists a (non-unique) function $\varepsilon \colon \ZZ/d \to \QQ$ such that, for all $\lambda \in (\ZZ/d)^\times$,
\begin{equation}\label{lem-del}
p(\lambda)=\sum_{a \in \ZZ/d} \varepsilon(a) \{\tfrac{a\lambda}{d} \}.
\end{equation}

\subsection{The Gross-Deligne conjecture} \

\noindent 
Let us keep the notation of the preceding section, with the additional assumption that $K=\alg$. Motivated by his geometric proof of the Chowla-Selberg formula, Gross made the following conjecture in \cite[p. 205]{Gro78}, the precise form of which was suggested by Deligne. 

\begin{conjecture}[Gross-Deligne]\label{conj-gd} For any function $\varepsilon \colon (\ZZ/d)^\times \to \QQ$ satisfying \eqref{lem-del}, the following relation holds: 
\begin{equation}\label{conjec}
P_\lambda(H) \sim_{\alg^\times} \prod_{a=1}^{d-1} \Gamma\left(1-\frac{a}{d} \right)^{\varepsilon(a/\lambda)}. 
\end{equation}
\end{conjecture}

\begin{remark}\label{koblitz-ogus} By a result of Koblitz and Ogus \cite{KO79}, the right hand sides of \eqref{conjec} corresponding to different choices of $\varepsilon$ agree up to an algebraic number, so it suffices to prove the conjecture for a single $\varepsilon$. Precisely, if a function $f \colon \ZZ/d \to \QQ$ satisfies $\sum_{a \in \ZZ/d} f(a)\{a\lambda \slash d\}=u$ for all $\lambda \in (\ZZ/d)^\times$, then the number
$$
\Gamma_{f, \lambda}:=(2\pi i)^{-u} \prod_{a=1}^{d-1} \Gamma\left(1-\frac{a}{d} \right)^{f(a \slash \lambda)}
$$ is algebraic for all $\lambda$. If, in addition, $f$ takes integer values, Gross and Koblitz show that $\Gamma_{f, \lambda}$ generates a Kummer extension of $\QQ(\xi)$ of degree dividing the least common multiple of $2$ and $d$ \cite[Thm. 4.5]{GK79}. 
\end{remark}

A strong evidence supporting the Gross-Deligne conjecture is that it holds in weight one, that is, for $H^1$ of abelian varieties with complex multiplication by an abelian number field, a theorem of Anderson \cite{And82}\footnote{Although this is not explicitly stated in \cite{And82}, it follows readily from Prop. 1.5 and Thm. 3.1 and 4.7 therein.}. On the contrary, the conjecture remains wide open in higher weight, even when the function $p$ is constant. In this case, the multiplication formula 
\begin{equation}\label{mult-gamma}
\prod_{a=1}^{d-1} \Gamma\left(1-\frac{a}{d} \right)=d^{-\frac{1}{2}}(2\pi)^{\frac{d-1}{2}}
\end{equation} allows to rewrite the statement as $P_\lambda(H) \sim_{\alg^\times} (2\pi i)^p$. Observe that this would follow from the Hodge conjecture or if all classes in $H$ were absolute Hodge. 
 
\subsection{Main results}\label{dethdg} \

\noindent
We now state the main results of the paper, which provide new evidence for the Gross-Deligne conjecture in higher weight. 

Let $d \geq 3$ be an integer, $\xi$ a primitive $d$-th root of unity and $K$ a subfield of $\alg$ containing $\QQ(\xi)$. Let $X$ be a smooth, projective variety, together with an automorphism $g$ of order $d$, with everything defined over $K$. 

We shall denote by $H^j(X, g)$ the largest sub-Hodge structure $V \subseteq H^j_B(X)$ such that $g^\ast$ acts on $V \otimes_\QQ \CC$ through primitive eigenvalues, that is
$$
H^j(X, g)=H^j_B(X) \otimes_{\QQ[\mu_d]} \QQ(\xi).
$$Then $H^j(X, g)$ has a natural structure of vector space over $\QQ(\xi)$, the dimension of which will be denoted by $h_j$. If $h_j \geq 1$, the $\QQ$-vector space $\det_{\QQ(\xi)} H^j(X, g)$ carries a pure Hodge structure of weight $j h_j$ with $\QQ(\xi)$-multiplication. 

Using the K\"unneth formula, $\det_{\QQ(\xi)} H^j(X, g)$ can be realized inside the cohomology of $X^{h_j}$, from which conditions \ref{h1} and \ref{h2} follow \cite[Lemma 2.1]{MR04}. Moreover, the associated function $p_j \colon (\ZZ/d)^\times \to \ZZ$ is given by
\begin{equation}\label{func-det}
p_j(\lambda)=\sum_{p+q=j} p h^{p, q}_\lambda(X),
\end{equation} where $h^{p, q}_\lambda(X)$ denotes the dimension of the $\xi^\lambda$-eigenspace of $H^{p, q}(X)$. 

We prove the following alternating variant of the Gross-Deligne conjecture\footnote{Strictly speaking, the period is only defined when $h_j \geq 1$. For notational convenience, we shall set $P_\lambda(\det_{\QQ(\xi)} H^j(X, g))=1$ if this is not the case.} for the Hodge structures $\mathrm{det}_{\QQ(\xi)} H^j(X, g):$ 

\begin{maintheorem}\label{maintheoA} There exists a function $\varepsilon: \ZZ/d \to \QQ$ satisfying the conditions: 
\begin{enumerate}
\item[(i)] The following equality holds for all $\lambda \in (\ZZ/d)^\times:$
$$
\sum_{a \in \ZZ/d} \varepsilon(a) \{\tfrac{a\lambda}{d} \}=\sum_{j=0}^{2n} (-1)^j \sum_{p+q=j} p h^{p, q}_\lambda(X). 
$$ 
\item[(ii)] For each $\lambda \in (\ZZ/d)^\times$, there exists an algebraic number $\alpha_\lambda$ such that 
$$
\prod_{j=0}^{2n} P_\lambda(\mathrm{det}_{\QQ(\xi)} H^j(X, g))^{(-1)^j} \sim_{L^\times} \alpha_\lambda \cdot \prod_{a=1}^{d-1} \Gamma\left(1-\frac{a}{d} \right)^{\varepsilon(a \slash \lambda)}, 
$$ where $L$ is the extension of $K$ generated by $i$ and $d^{\frac{1}{d-1}}$. Moreover, if $w$ denotes the least common multiple of $2$ and $d$, then $\alpha_\lambda^w \in K^\times$.\end{enumerate}
\end{maintheorem}

\begin{remark} The presence of the extension $L$ in the statement is only due to the extra factors needed to transform powers of $2\pi i$ into products of gamma values using the multiplication formula \eqref{mult-gamma}. In fact, the proof will show that there exist integers $M_\lambda$ such that, for all integer-valued function $\varepsilon \colon \ZZ/d \to \ZZ$ satisfying
$$
\sum_{a \in \ZZ/d} \varepsilon(a) \{\tfrac{a\lambda}{d}\}+M_\lambda=\sum_{j=0}^{2n} (-1)^j \sum_{p+q=j} p h^{p, q}_\lambda(X), 
$$ there exists algebraic numbers $\alpha_{\lambda, \varepsilon}$, with $\alpha_{\lambda, \varepsilon}^w \in K^\times$, such that: 
$$
\prod_{j=0}^{2n} P_\lambda(\mathrm{det}_{\QQ(\xi)} H^j(X, g))^{(-1)^j} \sim_{K^\times} \alpha_{\lambda, \varepsilon} (2\pi i)^{M_\lambda} \prod_{a=1}^{d-1} \Gamma\left(1-\frac{a}{d} \right)^{\varepsilon(a \slash \lambda)}.  
$$
\end{remark}

\begin{remark} Using techniques in Arakelov geometry, Maillot and R\"ossler prove that, for each odd primitive Dirichlet character $\chi: (\ZZ/d)^\times \to \CC$, the equality of complex numbers
\begin{align}\label{mai-ros}
\sum_{j=0}^{2n} (-1)^j \sum_{\lambda \in (\ZZ/d)^\times} &\log \abs{P_\lambda(\mathrm{det}_{\QQ(\xi)} H^j(X, g))} \chi(\lambda) \nonumber \\
&=\frac{L'(\chi, 0)}{L(\chi, 0)} \sum_{j=0}^{2n} (-1)^j \sum_{\lambda \in (\ZZ/d)^\times} \chi(\lambda) \sum_{p+q=j} p h^{p, q}_\lambda(X)
\end{align} holds, up to addition of a term of the form $\sum_{\lambda \in (\ZZ/d)^\times} \sum_i (b_{i, \lambda} \log \abs{\alpha_{i, \lambda}}) \chi(\lambda)$, where $\alpha_{i, \lambda} \in K^\times,$ $b_{i, \lambda} \in \QQ(\xi)$ and $i$ runs through a finite set of indices \cite[Thm. 1]{MR04}\footnote{But be aware that the square in Conj. $B(\Mca_0, N, \chi)$ in page 730 of \textit{loc. cit.} is a misprint.}, cf. also the Bourbaki seminar \cite{Sou07}. From \eqref{mai-ros} and a classical formula relating $L'(\chi, 0) \slash L(\chi, 0)$ to gamma values, the authors deduce a weak version of Theorem~\ref{maintheoA}, in which the period is replaced by its absolute value, under the additional assumption that $d$ is prime \cite[\S4]{MR04}. 
\end{remark}

Let us emphasize that separating the periods  coming from different degrees seems to be out of the reach of the methods used in this paper. Nevertheless, Theorem \ref{maintheoA} implies the original Gross-Deligne conjecture in some interesting cases where $H^j(X, g)$ vanishes in all degrees but one, as already observed by Maillot and R\"ossler \cite[Cor. 4.2]{MR04}.

\begin{corollary} Let $X$ be an $n$-dimensional smooth, complete intersection over $\alg$, together with an automorphism $g$ of finite order $d \geq 3$. The Gross-Deligne conjecture holds for the Hodge structure $\det_{\QQ(\xi)} H^n(X, g)$. 
\end{corollary}

\begin{proof} Let $\iota \colon X \hookrightarrow \PP^N$ be a projective embedding of $X$. By the Lefschetz hyperplane theorem, $\iota^\ast \colon H^j(\PP^N) \to H^j(X)$ is an isomorphism for all $j \neq n$. Since such a non-zero $H^j(X)$ is one-dimensional, $g^\ast$ acts by multiplication by a $d$-th root of unity, which can only be $\pm 1$ because integral cohomology is preserved. Taking into account that $d \geq 3$, the only non-trivial $H^j(X, g)$ is the one in degree $n$. 
\end{proof}

We now discuss two instances of Theorem \ref{maintheoA} which are easier to handle because the quotient of $X$ by the action of $g$ is a smooth variety, so one can directly adapt the approach to Euler's formula sketched in the introduction. In general, before taking the quotient, it is necessary to perform a sequence of equivariant blow-ups to reduce the computation to the case of a cyclic cover of a smooth variety ramified along a simple normal crossings divisor.  
\begin{proposition}\label{freeaction} If $g$ acts without fixed points on $X$, then
$$
\prod_{j=0}^{2n} P_\lambda(\mathrm{det}_{\QQ(\xi)} H^j(X, g))^{(-1)^j} \sim_{K^\times} \alpha_\lambda (2\pi i)^{n\chi(X) \slash 2d}. 
$$
\end{proposition}

We refer to Sect. \ref{proof-free} for a precise description of the constant $\alpha_\lambda$. Observe that, if we further assume that $H^j(X, g)$ vanishes in all degrees but one, which is then necessarily $j=n$ by Poincar\'e duality, then $\det_{\QQ(\xi)} H^n(X, g)$ consists entirely of Hodge classes (see Remark \ref{hodge-classes}). It is not known whether they are absolute Hodge, \textit{a fortiori} algebraic. Nonetheless, the above proposition shows that their periods are, up to an algebraic factor, the expected powers of $2\pi i$. 

In the case of curves, we recover a result due to Terasoma \cite[Thm. 2.1.3]{Ter92}. To formulate it, let $\pi: X \to Z$ denote the quotient by the action of $g$, $g(Z)$ the genus of $Z$, and $S \subset Z$ the branch locus of $\pi$. If $S=\coprod \mathrm{Spec}(K_i)$, we denote by $\delta_S$ the product of the square roots of the discriminants of the extensions $K_i \slash K$ (see Example \ref{sqrdisc}). We may choose a rational section $\omega$ of the sheaf of logarithmic differentials $\Omega^1_Z(\log S)$ with non-vanishing residues along $S,$ so that $\mathrm{div}(\omega)$ is supported on $Z \setminus S$. By Kummer theory, there exists a generator $f_\lambda$ of the cyclic field extension $K(X) \slash K(Z)$ corresponding to $\pi$ such that $f_\lambda^g=\xi^\lambda f_\lambda$. Then $F_\lambda=f_\lambda^d$ is a rational function on $Z$. Moreover, it is possible to find $f_\lambda$ invertible at the points over $\mathrm{div}(\omega)$. Let $z$ be a closed point in $Z \setminus S$ and $x$ a lifting to $X$. The value of $f_\lambda(x)$ modulo $K^\times$, when defined, depends only on $z$, so we shall write $f_\lambda(z)$. Consider the element in $\alg^\times \slash K^\times$ given by 
$$
\alpha_\lambda=\delta_S \cdot \prod_{z \in Z \setminus S} \ \prod_{\sigma \in \mathrm{Aut}_K(K(z))}f_\lambda(z^\sigma)^{\ord_z(\omega)} \cdot \prod_{s \in S} \Res_s(\omega)^{-\frac{1}{d} [K(s) \colon K] \ord_s(F_\lambda)}, 
$$ where the last term is defined as $\exp(-\frac{1}{d} [K(s) \colon K] \ord_s(F_\lambda) \log(\Res_s(\omega)))$ for some branch of the logarithm. Since $f_\lambda$ is unique up to multiplication by an element of $K(Z)^\times$ and $K$ contains $\QQ(\xi),$ the class $\alpha_\lambda$ is independent of all the above choices.  

\begin{maintheorem}\label{maintheoB} The following relation holds: 
$$
P_\lambda(\mathrm{det}_{\QQ(\xi)} H^1(X, g)) \sim_{K^\times} \alpha_\lambda \cdot (2\pi i)^{g(Z)-1} \cdot \prod_{s \in S} \Gamma\left(\frac{\ord_s(F_\lambda)}{d} \right)^{[K(s) \colon K]}. 
$$ 
\end{maintheorem}

\subsection{Analogue over finite fields} \

\noindent 
It is worth mentioning that Theorem \ref{maintheoB} admits an analogue over finite fields, in which epsilon factors play the role of periods and gamma values are replaced by certain Gauss sums. Let $\mathfrak{p}$ be a prime ideal of $K$ not dividing $d$, with residue field $\FF_q$. Assume that both $X$ and $g$ have good reduction at $\mathfrak{p}$. Then the Artin $L$-function $L(X_\mathfrak{p} \slash \FF_q, \lambda, T)$ satisfies a functional equation of the form
$$
L(X_{\mathfrak{p}} \slash \FF_q, \lambda, T)=\varepsilon(\lambda) T^{-\chi} L(X_{\mathfrak{p}} \slash \FF_q, -\lambda, (qT)^{-1}), 
$$ where $\chi$ is an integer and $\varepsilon(\lambda)$ a non-zero complex number. To give a formula for $\varepsilon(\lambda)$, let us denote by $\mathrm{Teich}_{\mathfrak{p}} \colon \mu_d(\FF_q) \to \mu_d(K)$ the inverse of the reduction map modulo $\mathfrak{p}$ and define, for $1 \leq a \leq d-1$, the Gauss sum
$$
\tau_{\mathfrak{p}}(a \slash d)=-\sum_{x \in \FF_q^\times} \mathrm{Teich}_{\mathfrak{p}}(x^{\frac{a}{d}(q-1)})\exp(\tfrac{2\pi i}{p} \cdot \Tr_{\FF_q \slash \FF_p}(x)). 
$$ 

\begin{maintheorem}\label{maintheoC} There exists a root of unity $\alpha_\mathfrak{p}$ of order $w$ such that
$$
\varepsilon(\lambda)=\alpha_{\mathfrak{p}} \cdot q^{g(Z)-1} \cdot \prod_{s \in S} \tau_\mathfrak{p}\left(\frac{\ord_s(F_\lambda)}{d}\right)^{[K(s) \colon K]}.  
$$ Moreover, if $F_{\mathfrak{p}} \in \Gal(\alg \slash K)^{\textup{ab}}$ denotes a geometric Frobenius at $\mathfrak{p}$, the following ``reciprocity law'' holds: $F_{\mathfrak{p}}(\alpha)=\alpha_\mathfrak{p} \alpha$. 
\end{maintheorem}

We conjecture that a similar result is true in any dimension. In a subsequent paper \cite{fres}, we shall use an $\ell$-adic avatar of our techniques to prove Theorem \ref{maintheoC}, as well as an alternating version of this conjecture. 

\section{The Saito-Terasoma theorem}\label{sec-2}

In this section, we recall the Saito-Terasoma theorem \cite{ST97}. The following notation will be used: $k$ and $F$ denote subfields of $\CC$, and $U$ a smooth quasi-projective variety, of dimension $n$, over $k$. By resolution of singularities, there exists a smooth projective variety $X$ over $k$ containing $U$ as the complement of a simple normal crossings divisor $D$. We fix such an $X$ throughout.

\subsection{The category $M_{k, F}(U)$}\label{sec-2.1} \

\noindent
Following \cite[Def. 2, p. 872]{ST97}, we consider the category of triples 
\begin{displaymath}
\Mca=((\Ec, \nabla), V, \rho),
\end{displaymath} where:
\begin{enumerate}
\item[$\bullet$] $\Ec$ is a locally free $\Oc_U$-module of finite rank endowed with an integrable connection $\nabla: \Ec \to \Ec \otimes_{\Oc_U} \Omega^1_U$ with regular singularities along $D$,
\smallskip
\item[$\bullet$] $V$ is a local system of $F$-vector spaces on $U_\CC^\textup{an}$,
\smallskip
\item[$\bullet$] $\rho: V \to \Ec^\textup{an}$ is a morphism of sheaves inducing an isomorphism of complex local systems $V \otimes_F \CC \simeq \ker(\nabla^\textup{an})$.
\end{enumerate}

Recall that $(\Ec, \nabla)$ has \textit{regular singularities along $D$} if there exists a locally free $\Oc_X$-sheaf $\Ec_X$ and a logarithmic connection $\nabla_X: \Ec_X \to \Ec_X \otimes_{\Oc_X} \Omega^1_X(\log D)$ satisfying $(\Ec_X, \nabla_X)_{| U} \simeq (\Ec, \nabla)$. Such a pair will be called a \textit{logarithmic extension}. We refer the reader to \cite{Del70} or \cite{Kat71} for basic properties of regular singularities, including the fact that Gauss-Manin connections are always regular singular. 

We shall denote by $\rk(\Mca)$ the rank of $\Ec$ and define the \textit{determinant} of $\Mca$ as the rank one object $\det(\Mca)=((\det_{\Oc_U} \Ec, \tr(\nabla)), \det_F, \det \rho)$. 

The category $M_{k, F}(U)$ possesses a ``unit object'' 
$\unit_U=((\Oc_U, d), F, \rho),$ where $\rho$ denotes the inclusion of the constant sheaf $F$ into $\Oc_{U}^\textup{an}$. 

The set $MPic_{k, F}(U)$ of isomorphism classes of rank one objects in $M_{k, F}(U)$ forms a group under tensor product, with the class of $\unit$ as unit. If $U=\mathrm{Spec}(L)$, for a finite extension $L$ of $k$, it has the following description: 
$$
MPic_{k, F}(L) \simeq k^\times \backslash (\CC^\times)^{\Hom_k(L, \CC)} \slash (F^\times)^{\Hom_k(L, \CC)}. 
$$

Indeed, rank one objects in $M_{k, F}(L)$ are triples $M=(M_L, M_F, (\rho_\sigma)_\sigma)$ consisting of one dimensional vector spaces $M_L$ over $L$ and $M_F$ over $F$, together with isomorphisms $\rho_\sigma \colon M_L \otimes_{L, \sigma} \CC \stackrel{\sim}{\to} M_F \otimes_F \CC$ for each $\sigma \in \Hom_k(L, \CC)$. Let us choose basis $m_L$ and $m_F$ and write $\rho_\sigma(m_L \otimes_{L, \sigma} 1_\CC)=\alpha_\sigma m_F \otimes_F 1_\CC$. Then the assignment $M \mapsto (\alpha_\sigma)_\sigma$ induces the desired isomorphism. 

Moreover, there is a norm map 
\begin{equation}\label{eq-norm}
N_{L \slash K} \colon MPic_{k, F}(L) \to MPic_{k, F}(k)
\end{equation} given, in the above notation, by $[M] \mapsto \prod_\sigma a_\sigma$. 

\subsection{Periods of connections}\label{sec-per-con} \

\noindent
An integrable connection $\nabla \colon \Ec \to \Ec \otimes_{\Oc_U} \Omega^1_U$ canonically extends to a complex of sheaves $(\Ec \otimes_{\Oc_U} \Omega^\bullet_U, \nabla)$ which will be denoted by $DR(\Ec, \nabla)$. By analytification, one deduces a complex of sheaves $DR(\Ec^\textup{an}, \nabla^\textup{an})$ on $U_\CC^\textup{an}$ and canonical maps of complex vector spaces 
\begin{equation}\label{eq-7}
\HH^\bullet(U, DR(\Ec, \nabla)) \otimes_k \CC \longrightarrow \HH^\bullet(U_\CC^\textup{an}, DR(\Ec^\textup{an}, \nabla^\textup{an})).
\end{equation}

When $(\Ec, \nabla)$ has regular singularities, this map is an isomorphism by a theorem of Deligne \cite[Chap. II, Thm. 6.2]{Del70}. Besides, $DR(\Ec^\textup{an}, \nabla^\textup{an})$ is a resolution of the local system $\ker(\nabla^\textup{an})$ \cite[Chap. I, Prop. 2.19]{Del70}, whence:
\begin{equation}\label{eq-8}
H^\bullet(U_\CC^\textup{an}, \ker(\nabla^\textup{an})) \stackrel{\sim}{\longrightarrow} \HH^\bullet(U_\CC^\textup{an}, DR(\Ec^\textup{an}, \nabla^\textup{an})). 
\end{equation}

If, in addition, $V$ and $\rho$ as in paragraph \ref{sec-2.1} are given, one gets isomorphisms 
\begin{equation*}
H^\bullet(\rho) \colon \ \HH^\bullet(U, DR(\Ec, \nabla)) \otimes_k \CC \stackrel{\sim}{\longrightarrow} H^\bullet(U_\CC^\textup{an}, V) \otimes_F \CC 
\end{equation*} from \eqref{eq-7}, \eqref{eq-8} and the map induced by $\rho$ on cohomology. For the unit object, one recovers the classical comparison between de Rham and Betti cohomology.  

Let $P_j$ denote the matrix of $H^j(\rho)$ with respect to any basis of $\HH^j(U, DR(\Ec, \nabla))$ over $k$ and any basis of $H^j(U_\CC^\textup{an}, V)$ over $F$. 

\begin{definition}
The \textit{determinant of periods} of $\Mca=((\Ec, \nabla), V, \rho)$ is given by the alternating product
$$
\per(\Mca)=\prod_{j=0}^{2n} (\det P_j)^{(-1)^j} \in k^\times \backslash \CC^\times \slash F^\times. 
$$
\end{definition}

\subsection{The Saito-Terasoma theorem}

\begin{theorem}[Saito-Terasoma]\label{st-theo} The following equality holds in $k^\times \backslash \CC^\times \slash F^\times:$
$$
\per(\Mca)=\per(\unit_U)^{\rk(\Mca)} \cdot (\det(\Mca), \cyc) \cdot \Gamma(\nabla^\vee \colon \Ec^\vee). 
$$
\end{theorem}

Here $(\Ec^\vee, \nabla^\vee)$ denotes the dual connection\footnote{Recall that $\nabla^\vee$ is uniquely determined by the requirement $d\langle e, e^\vee \rangle=\langle \nabla e, e^\vee \rangle + \langle e, \nabla^\vee e \rangle$, where $\langle \ , \ \rangle: \Ec \times \Ec^\vee \to \Oc_U$ denotes the canonical duality pairing.}. The gamma factor $\Gamma(\nabla^\vee \colon \Ec^\vee)$ and the pairing with the relative canonical cycle $\cyc$ will be introduced in paragraphs \ref{sec-gamma} and \ref{sec-pairing} below. We shall also give a formula for $\per(\unit_U)$. 

\subsection{The gamma factor}\label{sec-gamma} \

\noindent 
Let $(\Ec_X, \nabla_X)$ be a logarithmic extension of $(\Ec, \nabla)$ to $X$. For each irreducible component $D_i$ of $D$, let $k_i$ denote its field of constants, that is, the finite degree extension $\Gamma(D_i, \Oc_{D_i})$ of $k$. The composition
$$
\Ec_X \stackrel{\nabla_X}{\longrightarrow} \Ec_X \otimes \Omega^1_X(\log D) \xrightarrow{\id \otimes \Res_{D_i}} \Ec_X \otimes_{\Oc_X} \Oc_{D_i}=\colon \Ec_{X | D_i}
$$ induces an ($\Oc_{D_i}$-linear) endomorphism $\Res_{D_i} \nabla_X \in \End(\Ec_{X | D_i})$ called the \textit{residue of the connection along $D_i$}. Consider its characteristic polynomial
$$
\Phi_{\Ec_X, i}(T)=\det(T-\Res_{D_i} \nabla_X) \in k_i[T],
$$ the roots of which are classically called \textit{exponents} of $(\Ec_X, \nabla_X)$ at $D_i$. 

\begin{definition}
A logarithmic extension of $(\Ec, \nabla)$ is \textit{small} if none of the exponents at none of the components $D_i$ is an integer $n \leq 0$.
\end{definition}

\begin{remark}\label{small-ext}
Small extensions always exists: indeed, starting from any extension $(\Ec_X, \nabla_X)$, the sheaf $\Ec_X \otimes \Oc_X(\sum \mu_j D_j)$ carries a natural logarithmic connection with residue $\Res_{D_i} \nabla_X -\mu_i \cdot \id_{\Oc_{D_i}}$ along $D_i$ \cite[Lemma 2.7]{EV92}. 
\end{remark}

Given a polynomial $P \in k[T]$ such that $P(n) \neq 0$ for all integers $n \leq 0$, we put
$$
\Gamma(P)=\prod_{P(\alpha)=0} \Gamma(\alpha) 
$$ (roots counted with multiplicity). If the coefficients of $P$ lie instead in a finite extension $L$ of $k$, we shall first apply the norm $N_{L \slash k}: L(T) \to k(T)$.   

\begin{definition} The \textit{gamma factor of a small extension} $(\Ec_X, \nabla_X)$ of $(\Ec, \nabla)$ is the non-zero complex number\footnote{To eliminate any possible ambiguity, $ here \chi(D_i^\circ)=\sum_{j=0}^{2(n-1)} (-1)^j \dim_{k_i} H_{dR}^j(D_i^\circ \slash k_i)$.}
$$
\Gamma(\nabla_X \colon \Ec_X)=\prod_{i \in I} \Gamma(N_{k_i \slash k} \Phi_{\Ec_X, i})^{\chi(D_i^\circ)}. 
$$
\end{definition}

Observe that the exponents depend on the extension, whereas their classes modulo $\ZZ$ are uniquely determined by $(\Ec, \nabla)$. Indeed, if $\alpha$ is an exponent, then $\exp(-2\pi i \alpha)$ is an eigenvalue of the local monodromy of $\ker(\nabla^\textup{an})$ around $D_i$ \cite[Ch. II, Prop. 3.11]{Del70}. Together with the functional equation $\Gamma(s+1)=s\Gamma(s)$, this implies that the class of $\Gamma(\nabla_X \colon \Ec_X)$ modulo $k^\times$ does not depend on the choice of the small extension \cite[Lemma 1.8]{ST97}.

\begin{definition} The \textit{gamma factor} of $(\Ec, \nabla)$ is the class of $\Gamma(\nabla_X \colon \Ec_X)$ modulo $k^\times$ for any small extension $(\Ec_X, \nabla_X)$. It will be denoted by\footnote{Saito and Terasoma use the notation $\Gamma(\nabla \colon \Mca)$; we prefer to emphasize that the gamma factor depends only on $(\Ec, \nabla)$.} $\Gamma(\nabla \colon \Ec)$. 
\end{definition}

\begin{remark}\label{gam-brie} When $k$ and $F$ are subfields of $\alg$, Brieskorn \cite[Satz 4]{Bri70} shows that, if $\alpha$ is an exponent of $(\Ec_X, \nabla_X)$, then both $\alpha$ and $\exp(-2\pi i \alpha)$ are algebraic numbers. Hence, by the Gelfond-Schneider theorem, all $\alpha$ are rational (cf. also \cite[Lemma 3.1]{Sim90}). In this case, the gamma factor takes the form
$$
\Gamma(\nabla \colon \Ec)=\prod_{i \in I} \prod_{j=1}^{\rk(\Ec)} \Gamma(\alpha_{i j})^{[k_i \colon k] \chi(D_i^\circ)}. 
$$ 
\end{remark}

\subsection{The pairing with the relative canonical cycle}\label{sec-pairing} \

\noindent
In this section, we introduce the remaining terms in the Saito-Terasoma theorem. Recall that $n$ always denotes the dimension of $X$.

\subsubsection{The relative Chow group of zero-cycles} \

\noindent 
Let us denote by $\mathcal{K}_{n, X}$ (resp. $\mathcal{K}_{n, D_i}$) the Zariski sheaf on $X$ (resp. on $D_i$) associated to Quillen's $K$-theory in degree $n$ \cite[\S 3]{Qui73}, and consider the complex 
$$
\mathcal{K}_n(X \ \mathrm{mod} \ D)=[\mathcal{K}_{n, X} \longrightarrow \bigoplus_{i \in I} (\iota_i)_\ast \mathcal{K}_{n, D_i}], 
$$ where $\mathcal{K}_{n, X}$ sits in degree zero and $\iota_i \colon D_i \hookrightarrow X$ is the inclusion\footnote{Note that the definition of $\mathcal{K}_n(X \ \mathrm{mod} \ D)$ makes sense for any family $D=(D_i)_{i \in I}$ of regular closed subschemes of $X$, not necessarily forming a simple normal crossings divisor.}. 

\begin{definition}
The \textit{relative Chow group of zero cycles} is the hypercohomology of this complex: 
$$
\CH^n(X \ \mathrm{mod} \ D)=\HH^n(X, \mathcal{K}_n(X \ \mathrm{mod} \ D)).
$$
\end{definition}

By the Bloch-Quillen formula $\CH^p(X) \simeq H^p(X, \mathcal{K}_{p, X})$ \cite[\S7, Thm. 5.19]{Qui73}, this group fits into a long exact sequence 
\begin{equation}\label{blo-qui}
\cdots \longrightarrow \bigoplus_{i \in I} H^{n-1}(D_i, \mathcal{K}_{n, D_i}) \longrightarrow \CH^n(X \ \mathrm{mod} \ D) \longrightarrow \CH^n(X) \longrightarrow 0. 
\end{equation}

Saito and Terasoma prove that $\CH^n(X \ \mathrm{mod} \ D)$ admits an adelic presentation
\begin{equation}
\CH^n(X \ \mathrm{mod} \ D)=\coker\left(\partial \colon \bigoplus_{y \in X_1} B_y \longrightarrow \bigoplus_{x \in X_0} A_x \right), 
\end{equation} where $X_i$ denotes the set of points of $X$ of dimension $i$ \cite[Prop. 1]{ST97}. In this formula, the group $A_x$ is an extension of $\ZZ$ by $\bigoplus_{i \in I_x} k(x)^\times$, with $I_x \subseteq I$ the index subset of irreducible components of $D$ such that $x \in D_i$; the group $B_y$ is an extension of $k(y)^\times$ by $\bigoplus_{i \in I_y} K_2(k(y))$, with $I_y$ similarly defined; and the morphism $\partial$ is the direct sum of maps $\partial_{x, y} \colon B_y \to A_x$, for $x$ in the closure of $y$, induced from the order $\ord_x \colon k(y)^\times \to \ZZ$ and the tame symbol $(\ , \ )_x \colon K_2(k(y)) \to k(x)^\times$. 

In particular, each closed point of $U$ defines a class $[x]$ in $\CH^n(X \ \mathrm{mod} \ D)$, the image of $1 \in \ZZ=A_x$. A corollary of the adelic presentation is that the relative Chow group of zero cycles is actually generated by these classes, as $x$ runs through any dense open subscheme of $U$ \cite[Cor. p. 880]{ST97}. 

\subsubsection{The relative canonical cycle} 

\begin{definition} A \textit{partial trivialization} of a locally free $\Oc_X$-module $\Fca$ along $D$ is a family $r=(r_i)_{i \in I}$ of surjective morphisms $r_i \colon \Fca_{| D_i} \to \Oc_{D_i}$. 
\end{definition}

For example, the Poincar\'e residues $\mathrm{Res}_{D_i} \colon \Omega^1_X(\log D)_{| D_i} \to \Oc_{D_i}$ provide a partial trivialization of the sheaf of logarithmic differentials $\Omega^1_X(\log D)$. 

Let $(\Fca, r)$ be a locally free $\Oc_X$-module of rank $n$, together with a partial trivialization $r$ along $D$. By a construction due to Saito \cite[\S1]{Sai93}, to such a pair it is attached a relative top Chern class $c_n(\Fca, r) \in \CH^n(X \ \mathrm{mod} \ D)$ mapping to the usual top Chern Class under \eqref{blo-qui}. Let us briefly recall how it is defined. 

Consider the geometric vector bundle $V=\mathbf{V}(\Fca)$ associated to $\Fca$ and the morphisms $r_i \colon V_{D_i} \to \mathbb{A}^1_{D_i}$ induced by the partial trivialization. Define $\Delta_i=r_i^{-1}(\{1\})$, where $\{1\} \subset \mathbb{A}^1_{D_i}$ denotes the unit section. Then $\Delta=(\Delta_i)_{i \in I}$ is a family of regular closed subschemes of $V$ disjoint from the zero section $\{0\} \subset V$. Thus the local cohomologies $H^n_{\{0\}}(V, \mathcal{K}_{n, V})$ and $\HH^n_{\{0\}}(V, \mathcal{K}_n(V \ \mathrm{mod} \ \Delta))$ are isomorphic. Besides, since $V \to X$ and $\Delta_i \to D_i$ are affine bundles, there is an isomorphism 
$$
\HH^n(V, \mathcal{K}_n(V \ \mathrm{mod} \ \Delta)) \simeq \HH^n(X, \mathcal{K}_n(X \ \mathrm{mod} \ D))
$$ by the homotopy property of $K$-theory \cite[\S7, Prop. 4.1]{Qui73}. This and the natural map $\HH^n_{\{0\}}(V, \mathcal{K}_n(V \ \mathrm{mod} \ \Delta)) \to \HH^n(V, \mathcal{K}_n(V \ \mathrm{mod} \ \Delta))$ induce a morphism
\begin{equation}\label{kt}
H^n_{\{0\}}(V, \mathcal{K}_{n, V}) \longrightarrow \CH^n(X \ \mathrm{mod} \ D).  
\end{equation} 

With this notation, the \textit{relative top Chern class} $c_n(\Fca, r)$ is defined as the image of the cycle class $\{0\} \in H^n_{\{0\}}(V, \mathcal{K}_{n, V})$ under \eqref{kt}. Specializing to $\Omega^1_Z(\log D),$ together with the partial trivialization $\Res=(\Res_{D_i})$, we get: 

\begin{definition} The \textit{relative canonical cycle} is 
$$
\cyc=(-1)^n c_n(\Omega^1_X(\log D), \Res) \in \CH^n(X \ \mathrm{mod} \ D). 
$$
\end{definition}

\subsubsection{The pairing}\

\noindent

Let $x$ be a closed point of $U$ and $\Mca=((\Ec, \nabla), V, \rho)$ a rank one object in $M_{k, F}(U)$. The \textit{fiber of $\Mca$ at $x$} is the rank one object 
$$
\Mca_{| x}=(\Ec_x, V_x, (\rho_{x, \sigma})_{\sigma \in \Hom_k(k(x), \CC)}) \in M_{k, F}(k(x)), 
$$ where $\Ec_x$ and $V_x$ denote the fibers at $x$, and $\rho_{x, \sigma} \colon V_x \to \Ec_x^{\textup{an}}$ is obtained by restriction of $\rho$. Taking the norm \eqref{eq-norm} down to $k$, we get an element 
$$
N_{k(x) \slash k}(\Mca_{| x}) \in k^\times \backslash \CC^\times \slash F^\times. 
$$

\begin{proposition}[Saito-Terasoma]
There exists a unique pairing 
$$
(\ , \ ) \colon MPic_{k, F}(U) \times \CH^n(X \ \mathrm{mod} \ D) \longrightarrow k^\times \backslash \CC^\times \slash F^\times 
$$ such that $(\Mca, [x])=N_{k(x) \slash k}(\Mca_{| x})$ for all closed points $x$ of $U$. 
\end{proposition}

The uniqueness of the pairing follows from the fact that $\CH^n(X \ \mathrm{mod} \ D)$ is generated by the classes $[x]$; the existence is the main result of Sects. 2 and 3 of \cite{ST97}. All the factors involved in the Saito-Terasoma theorem are now defined. 

In concrete terms, we may compute $N_{k(x) \slash k}(\Mca_{| x})$ as follows \cite[p. 882]{ST97}: choose a local basis $e$ of $\Ec$ around $x$ and a local basis $v$ of $V$ on a sufficiently small neighborhood $W$ of $x$ in $U_\CC^\textup{an}$. Then the function $\varphi=\rho(v) \slash e \colon W \to \CC$ is invertible and holomorphic, and we have
\begin{equation}\label{explicit-pairing}
N_{k(x) \slash k}(\Mca_{| x})=\prod_{\sigma \in \Hom_k(k(x), \CC)} \varphi(\sigma(x))^{-1}. 
\end{equation}

\begin{example}\label{pairing-curves} Let us make, after \cite[p. 392]{Sai93} and \cite[p. 882]{ST97}, all the constructions explicit in the case where $X$ is a curve and $D$ a reduced divisor. Then the relative Chow group of zero-cycles $\CH^1(X \ \mathrm{mod} \ D)$ coincides with the usual divisor class group with modulus
\begin{equation}\label{modulus-curve}
\CH^1(X \ \mathrm{mod} \ D) \simeq \coker(k(X) \longrightarrow \bigoplus_{x \in U_0} \ZZ \oplus \bigoplus_{x \in D} k(X)^\times \slash (1+\mathfrak{m}_x)).  
\end{equation} 

To compute the relative canonical cycle, we may choose a rational section $\omega$ of $\Omega^1_X(\log D)$ with non-vanishing residues along\footnote{If $\deg D\geq 2$, we can even choose a global section.} $D$, so the divisor of $\omega$ does not meet $D$. Then the image of $\cyc$ under the isomorphism \eqref{modulus-curve} is
$$
\cyc=-\mathrm{div}(\omega)-\sum_{x \in D} [\mathrm{Res}_x(\omega)], 
$$ where $[\mathrm{Res}_x(\omega)]$ denotes the class of $\Res_x(\omega) \in k^\times$ in $k(X)^\times \slash (1+\mathfrak{m}_x).$  

Now assume that a rank one object $\Mca$ in $M_{k, F}(U)$ is given, and let $\nabla_X$ denote a small extension of the underlying connection. Then: 
$$
(\Mca, \cyc)=\prod_{x \in U_0} N_{k(x) \slash k}(\Mca_{| x})^{-\ord_x(\omega)} \cdot \prod_{x \in D} \Res_x(\omega)^{\Tr_{k(x) \slash k} (\Res_x \nabla_X)}, 
$$ where the last term is defined as $\exp(\Tr_{k(x) \slash k}(\Res_x \nabla_x) \cdot \log(\mathrm{Res}_x(\omega)))$ for some branch of the logarithm. It is easy to see that its class in $k^\times \backslash \CC^\times \slash F^\times$ is independent of the choice of both the extension and the branch of the logarithm. 
\end{example}

\subsection{Variant with compact support}\label{stconsup} \

\noindent
The Saito-Terasoma theorem admits the following variant. If $(E, \nabla)$ is an integrable connexion on $U$ with regular singularities along $D$, we define its \textit{cohomology with compact support} as
$$
H_c^\bullet(U, DR(\Ec, \nabla))=\mathbb{H}^\bullet(X, DR(\Ec_X, \nabla_X))
$$ for any small extension $(\Ec_X, \nabla_X)$ of $(\Ec, \nabla)$ to $X$ \cite[Lemma 1.2]{ST97}. The name is justified by the existence of a canonical isomorphism 
$$
H^\bullet_c(U, DR(\Ec, \nabla)) \otimes_k \CC \stackrel{\sim}{\longrightarrow} H_c^\bullet(U_\CC^\textup{an}, \ker(\nabla^\textup{an})). 
$$

As in paragraph \ref{sec-per-con}, for each object $\Mca$ of $M_{k, F}(U)$ we deduce isomorphisms 
$$
H_c^\bullet(\rho) \colon \ H^\bullet_c(U, DR(\Ec, \nabla)) \otimes_k \CC \stackrel{\sim}{\longrightarrow} H_c^\bullet(U_{\CC}^\textup{an}, V) \otimes_F \CC.  
$$ The period $\per_c(\Mca) \in k^\times \backslash \CC^\times \slash F^\times$ is then defined as the alternating determinant of $H^\bullet_c(\rho)$ with respect to rational bases. From Theorem \ref{st-theo} for the dual object and Poincar\'e duality \cite[Lemma 1.7]{ST97}, we immediately derive:  

\begin{theorem}[Saito-Terasoma]\label{st-compsup} The following equality holds in $k^\times \backslash \CC^\times \slash F^\times:$
$$
\per_c(\Mca)=\per_c(\unit_U)^{\rk(\Mca)} \cdot  (\det(\Mca), \cyc) \cdot \Gamma(\nabla \colon \Ec)^{-1}.
$$
\end{theorem}

\subsection{Periods of the unit object} \

\noindent
In this final section, we give a formula for the determinant of periods of the unit object $\unit_U$. To formulate the result, let $D(m)$ denote the disjoint union of all schematic intersections $D_J=\bigcap_{i \in J} D_i$ indexed by subsets $J \subseteq I$ of cardinal $m$, with the usual convention $D(0)=X$. Observe that $D(m)$ is a smooth projective variety of dimension $n-m$. 

\begin{proposition}\label{prop-unit} The following equality holds in $k^\times \backslash \CC^\times \slash F^\times:$ 
\begin{align*}
\per(\unit_U)^2=(2\pi i)^{\sum_{m \geq 0} (-1)^m (n+m) \chi(D(m))}. 
\end{align*}
\end{proposition}

As a consequence of Poincar\'e duality: 

\begin{corollary}\label{prop-unit-cs} The following equality holds in $k^\times \backslash \CC^\times \slash F^\times:$
$$
\per_c(\unit_U)^2=(2\pi i)^{\sum_{m \geq 0} (-1)^m (n-m) \chi(D(m))}. 
$$
\end{corollary}

To prove the proposition, we can plainly assume that $F=\QQ$. Using a standard technique in mixed Hodge theory, we shall reduce to the smooth projective case. Recall that any smooth variety $Y$ over $k$ provides objects 
$$
H^j(Y)=(H^j_{dR}(Y \slash K), H^j_B(Y), \rho^j)
$$ of the abelian category $M_{k, \QQ}(k)$. Given an integer $r \geq 1$, we define the Tate twist $H^j(Y)(-r)$ as the tensor product $H^j(Y) \otimes H^2(\PP^1)^{\otimes r}$. As usual, we have
$$
\per(H^j(Y)(-r))=(2\pi i)^{r \dim H^j(Y)} \cdot \per(H^j(Y)).
$$

With the above notation, the identity 
$$
\sum_{j=0}^{2n} (-1)^j [H^j(U)]=\sum_{m \geq 0} (-1)^m \sum_{\ell=0}^{2(n-m)} (-1)^\ell [H^\ell(D(m))(-m)]
$$ holds in the Grothendieck group $K_0(M_{k, \QQ}(k))$, see \cite[p. 99]{PS08} for the corresponding statement at the Betti level, which generalizes at once. Since the determinant is multiplicative on exact sequences, it follows that 
\begin{equation}\label{k0mk}
\per(\unit_U)=(2\pi i)^{\sum_{m \geq 0} (-1)^m m \chi(D(m))} \cdot \prod_{m \geq 0} \per(\unit_{D(m)})^{(-1)^m}. 
\end{equation}

Thus we are reduced to prove the following lemma: 

\begin{lemma}\label{lem-per} Let $X$ be an $n$-dimensional smooth projective variety over $k$. Then:
$$
\per(\unit_X)^2=(2\pi i)^{n \chi(X)}. 
$$
\end{lemma}

\begin{proof} We may assume without loss of generality that $X$ is geometrically connected. Let us choose a hyperplane section $H$ of $X$ and consider its cohomology classes $\cl_{dR}(H) \in H^2_{dR}(X \slash k)$ and $\cl_B(H) \in H^2_B(X)$. Then $H^{2n}(X)$ is generated by $\cl(H)^n$. Since the comparison isomorphism sends $\cl_{dR}(H)$ to $(2\pi i) \cl_B(H)$ \cite[pp. 20-21]{Del82}, one has $[\det H^n(X)]=(2\pi i)^n$. Now, for each $0 \leq j \leq n$, the cup product $H^j(X) \otimes H^{2n-j}(X) \to H^{2n}(X)$ is non-degenerate (Poincar\'e duality) and compatible with the comparison isomorphism. It follows that
$$
[\det H^j(X)] \cdot [\det H^{2n-j}(X)]=(2\pi i)^{n \dim H^j(X)}, 
$$ and regrouping terms yields $\per(\unit_X)^2=(2\pi i)^{n\chi(X)}$. 
\end{proof}

\begin{remark}
A closer look at the proof of Lemma \ref{lem-per} gives in fact the identity
\begin{equation}\label{drB}
\per(\unit_X)=\sqrt{\disc(q_{dR}) \cdot \disc(q_B)} \cdot (2\pi i)^{\frac{n \chi(X)}{2}}, 
\end{equation} where $q_{dR}$ and $q_B$ denote the cup-product on $H^n_{dR}(X \slash k)$ and $H^n_B(X)$ respectively. For odd $n$, these bilinear forms are alternating, so the discriminant is a square and $\per(\unit_X)=(2\pi i)^{\frac{n\chi(X)}{2}}$. Observe that $n\chi(X) \slash 2$ is always an integer. 
\end{remark}

\begin{example}\label{sqrdisc} Let $L$ be an extension of $k$ of degree $d$, and let $\delta_{L \slash k}$ denote the determinant of the matrix $(\sigma_i(a_j))_{1 \leq i, j \leq d}$, where $a_1, \ldots, a_d$ is a basis of $L$ over $k$ and $\sigma_1, \ldots, \sigma_d$ are the complex embeddings of $L$ inducing the identity of $k$. The class of $\delta_{L \slash k}$ modulo $k^\times$ is independent of the chosen basis and $\delta_{L \slash k}^2$ agrees with the discriminant of $L$ over $k$. If $X$ denotes $\mathrm{Spec}(L)$, viewed as a zero-dimensional variety over $k$, then\footnote{Indeed, $H^0_{dR}(L \slash k)=L$ and $H_B(L)=\QQ^{\Hom_k(L)}$, so $(\sigma_i(a_j))_{1 \leq i, j \leq n}$ is the matrix of the comparison isomorphism with respect to these bases. Note that this is compatible with \eqref{drB} since $q_{dR}$ and $q_B$ are, respectively, the forms $\mathrm{Tr}(xy)$ and $\sum x_i^2$. } $\per(\unit_X)=\delta_{L \slash k}$. 
\end{example}

\begin{example}\label{sqrcurves} Let us now assume that $X$ is a curve of genus $g(X)$, write the ``divisor at infinity'' as $D=\coprod_{i \in I} \mathrm{Spec}(k_i)$ and define $\delta_D$ as the product of all $\delta_{k_i \slash k}$. Then it follows from equation \eqref{k0mk} and the previous example that  
$$
\per(\unit_U)=\delta_D (2\pi i)^{1-g(X)-\sum_{i \in I} [k_i \colon k]}. 
$$ Observe that the extension $k(\delta_D) \slash k$ corresponds to the kernel of the signature of the permutation action of $\Gal(\bar{k} \slash k)$ on $D(\bar{k})$. 
\end{example}

\section{Periods of cyclic covers}\label{sec-3}

In this section, we apply the Saito-Terasoma theorem to compute the periods of certain cyclic covers introduced by Esnault and Viehweg. From this we derive Proposition \ref{freeaction} and Theorem \ref{maintheoB}. Along the way, we shall prove some new cases of the Gross-Deligne conjecture (see Corollary \ref{gdforcovers} below). 

\subsection{Cyclic covers \`a la Esnault-Viehweg}\label{cyc-EV} \

\noindent 
The main reference for this section is \cite[\S 3]{EV92}. Let us assume that we are given the following data: 

\begin{itemize}
\item $d \geq 2$ is an integer, $\xi$ is a primitive $d$-th root of unity and $K$ is a subfield of $\alg$ containing $F=\QQ(\xi)$. 

\item $Z$ is a smooth projective variety, of dimension $n$, over $K$. 

\item $D=\sum_{i \in I} a_i D_i$ is an effective divisor with simple normal crossings on $Z$. 

\item $\Lcal$ is an invertible sheaf on $Z$ such that $\Oc_Z(D) \simeq \Lcal^d$, and $s \in H^0(Z, \Lcal^d)$ is a section with zero divisor $D$. 

\end{itemize}

We shall denote by $U$ the complement of $D_{\mathrm{red}}$ in $Z$ and, for each $0 \leq \lambda \leq d-1$, define the reduced divisor $D^{(\lambda)}=\sum_{a_i \lambda \nmid d} D_i$ and its complement $U^{(\lambda)}=Z\setminus D^{(\lambda)}$. 

Recall from \cite[\S3.5, 3.14]{EV92} that, out of these data, one can form a $d$-fold cyclic cover $\pi \colon Y \to Z$ ramified at $D$. A possible construction goes as follows: consider the geometric line bundles corresponding to $\Lcal^{-1}$ and $\Lcal^{-d}$, and let $\tau: \mathbf{V}(\Lcal^{-1}) \to \mathbf{V}(\Lcal^{-d})$ denote the ``$d$-th power map''. Then $s$ induces a geometric section $\sigma: Z \to \mathbf{V}(\Lcal^{-d})$. The associated cover $\pi \colon Y \to Z$ is the normalization of $\tau^{-1}(\sigma(Z))$, together with the restriction of the projection $\mathbf{V}(\Lcal^{-1}) \to Z$. 

A local computation shows that $Y$ has at most finite quotient singularities, lying above the singular locus of $D$ \cite[\S3.15, 3.24]{EV92}. Besides, $\pi \colon Y \to Z$ is a Galois cover, with cyclic Galois group $G$ of order $d$. Let us choose a generator $g$ and identify the group of characters with $\ZZ/d$ via the root of unity $\xi$. 

\begin{remark} Although $Y$ is not smooth, the construction of $\det_{\QQ(\xi)} H^j(Y, g)$ can be still carried out as in paragraph \ref{dethdg}. Indeed, since $Y$ has quotient singularities, any desingularization $h \colon \widetilde{Y} \to Y$ induces an injective morphism of Hodge structures $h^\ast \colon H^j_B(Y) \hookrightarrow H^j_B(\widetilde{Y})$ \cite[Thm. 8.2.4]{Del74} and the comparison isomorphism still holds, provided that de Rham cohomology is defined as the hypercohomology of the complex $\widetilde{\Omega}^\bullet_Y:=h_\ast \Omega^\bullet_{\widetilde{Y}}$ \cite[\S1]{Ste77}.   
\end{remark}

We shall compute the periods of $\det_{\QQ(\xi)} H^j(Y, g)$ by means of rank one connections on $Z$. For each $0 \leq \lambda \leq d-1$, define the invertible sheaf 
\begin{equation}
\Lcal^{(\lambda)}=\Lcal^{\lambda} \otimes_{\Oc_Z} \Oc_Z(-\sum_{i \in I} [\tfrac{a_i \lambda}{d}] D_i). 
\end{equation} In particular, $\Lcal^{(0)}=\Oc_Z$ and $\Lcal^{(\lambda)} \otimes_{\Oc_Z} \Lcal^{(d-\lambda)} \simeq \Oc_Z(D^{(\lambda)})$ for $\lambda \neq 0$. 

The following theorem summarizes the results that we need from \cite[\S 3]{EV92}: 

\begin{theorem}[Esnault-Viehweg]\label{esn-vie}
\begin{enumerate}
\item[] 

\item[\namedlabel{EVa}{(a)}] The sheaf $\Lcal^{(\lambda)^{-1}}$ is isomorphic to the subsheaf of $\pi_\ast \Oc_Y$ on which $g$ acts by multiplication by $\xi^\lambda$. Therefore,
$$
\pi_\ast \Oc_Y=\bigoplus_{\lambda=0}^{d-1} \Lcal^{(\lambda)^{-1}}. 
$$

\item[\namedlabel{EVb}{(b)}] Similarly, the eigenspace decomposition of $\pi_\ast \widetilde{\Omega}^p_Y$ is given by
$$
\pi_\ast \widetilde{\Omega}^p_Y=\bigoplus_{\lambda=0}^{d-1} \Lcal^{(\lambda)^{-1}} \otimes_{\Oc_Z} \Omega^p_Z(\log D^{(\lambda)}). 
$$ In particular, $h^{p, q}_\lambda(Y)=h^q(Z, \Lcal^{(\lambda)^{-1}} \otimes_{\Oc_Z} \Omega^p_Z(\log D^{(\lambda)})).$  

\item[(c)] The direct image of the exterior differential $\mathcal{O}_Y \to \widetilde{\Omega}^1_Y$ decomposes, according to (a) and (b), as a direct sum of logarithmic integrable connections
$$
\nabla^{(\lambda)} \colon \Lcal^{(\lambda)^{-1}} \longrightarrow \Lcal^{(\lambda)^{-1}} \otimes_{\Oc_Z} \Omega^1_Z(\log D^{(\lambda)})
$$ with residue $\mathrm{Res}_{D_i} \nabla^{(\lambda)}=\{ \frac{a_i \lambda}{d} \} \cdot \id_{\Oc_{D_i}}$ along $D_i$. 
\end{enumerate}
\end{theorem}

Since $\nabla^{(\lambda)}$ is a logarithmic connection on the smooth projective variety $Z$, with poles only along $D^{(\lambda)}$, its restriction to the complement $U^{(\lambda)}$ has regular singularities\footnote{This is a special case, for a finite morphism, of the aforementioned result that Gauss-Manin connections are regular singular.}. To get an object of the category $M_{K, F}(U^{(\lambda)})$, we observe that the direct image of the constant sheaf $F$ (resp. $\CC$) on $Y_\CC^\textup{an}$ also carries an action of $G$. For each $\lambda$, let $(\pi_\ast F)_\lambda$ (resp. $(\pi_\ast \CC)_\lambda$) denote the subsheaf of $\pi_\ast F$ (resp. $\pi_\ast \CC$) on which $g$ acts by multiplication by $\xi^\lambda$; then $(\pi_\ast F)_\lambda$ restricts to a rank one local system on $U^{(\lambda)}$. The analytic horizontal sections of $\Lcal^{(\lambda)^{-1}}$ are canonically identified with $(\pi_\ast \CC)_\lambda$, and the inclusion $\rho_\lambda \colon (\pi_\ast F)_\lambda \hookrightarrow (\pi_\ast \CC)_\lambda$ induces an isomorphism $(\pi_\ast F)_\lambda \otimes_F \CC \simeq \ker(\nabla^{(\lambda), \textup{an}}).$ Therefore, 
\begin{equation}\label{emlam}
\Mca_\lambda=((\Lcal^{(\lambda)^{-1}}, \nabla^{(\lambda)}), (\pi_\ast F)_\lambda, \rho_\lambda)_{| U^{(\lambda)}}
\end{equation} is a well defined object of $M_{K, F}(U^{(\lambda)})$. Taking into account Theorem \ref{esn-vie} and the description of the Hodge structure $\det_{\QQ(\xi)} H^j(Y, g)$, we have: 
\begin{equation*}
\per(\Mca_\lambda)=\prod_{j=0}^{2n} P_\lambda(\mathrm{det}_{\QQ(\xi)} H^j(Y, g))^{(-1)^j}. 
\end{equation*}

\subsection{Periods of cyclic covers} \

\noindent 
The main result of this section is the following: 

\begin{maintheorem}\label{maintheoD} For each $\lambda$, define an integer 
$$
M_\lambda=\frac{1}{2} \sum_{m \geq 0} (-1)^m (n+m) \chi(D^{(\lambda)}(m)).
$$ 
Let $\varepsilon \colon \ZZ/d \to \ZZ$ be any function satisfying, for all $\lambda \in (\ZZ/d)^\times,$  
\begin{equation}\label{rel-i}
\sum_{a \in \ZZ/d} \varepsilon(a) \{\tfrac{a \lambda}{d} \}+M_\lambda=\sum_{j=0}^{2n} (-1)^j \sum_{p+q=j} p h^{p, q}_\lambda(Y). 
\end{equation} Then, for each $\lambda \in (\ZZ/d)^\times$, there exists an algebraic number $\alpha_{\lambda, \varepsilon}$ such that 
\begin{equation}\label{rel-ii}
\prod_{j=0}^{2n} P_\lambda(\mathrm{det}_{\QQ(\xi)} H^j(Y, g))^{(-1)^j} \sim_{K^\times} \alpha_{\lambda, \varepsilon} \cdot (2\pi i)^{M_\lambda} \cdot \prod_{a=1}^{d-1} \Gamma\left(1-\frac{a}{d} \right)^{\varepsilon(a/\lambda)}. 
\end{equation} Moreover, $\alpha_{\lambda, \varepsilon}^w \in K^\times$, where $w$ is the least common multiple of $2$ and $d$.
\end{maintheorem}

\begin{corollary}\label{gdforcovers} Assume that the sheaf $\Lcal$ is ample and that the multiplicities $a_i$ are coprime to $d$. Then the Gross-Deligne conjecture holds for $\det_{\QQ(\xi)} H^n(Y, g).$ 
\end{corollary}

\begin{proof} By \cite[Lemma 1.5]{EV86} or \cite[Cor. 1.6]{Ara14}, these assumptions imply the vanishing of $H^n(Y, g)$ for all $j \neq n$. The statement is then a direct consequence of Theorem \ref{maintheoD} and the multiplication formula \eqref{mult-gamma}. 
\end{proof}

\subsection{Proof of Theorem \ref{maintheoD}} \

\noindent 
Thanks to Remark \ref{koblitz-ogus}, it suffices to prove the statement for only one function. We proceed in two steps. We first apply the Saito-Terasoma theorem to the object $\Mca_\lambda$ introduced in \eqref{emlam}; this yields the relation \eqref{rel-ii} for a certain function $\varepsilon_0$ to be described in the course of the proof. In a second step, we show that the equality \eqref{rel-i} is satisfied for this choice of $\varepsilon_0$. 

\subsubsection*{Step 1} Since the connections $\nabla^{(\lambda)}$ and $\nabla^{(d-\lambda)}$ are dual to each other over $U^{(\lambda)}$, the Saito-Terasoma theorem gives the following identity in $\CC^\times \slash K^\times$: 
$$
\per(\Mca_{\lambda})=\per(\unit_{U^{(\lambda)}}) \cdot (\Mca_{\lambda}, c_{Z \mathrm{mod} D^{(\lambda)}}) \cdot \Gamma(\nabla^{(d-\lambda)} \colon \Lcal^{(d-\lambda)^{-1}}). 
$$ 

By Proposition \ref{prop-unit}, there exists an algebraic number $C$, with $C^2 \in K^\times$, such that $\per(\unit_{U^{(\lambda)}})=C (2\pi i)^{M_\lambda}$. Let us define the number $\alpha_\lambda=C \cdot (\Mca_{\lambda}, c_{Z \mathrm{mod} D^{(\lambda)}})$ and the function $\varepsilon_0 \colon \ZZ/d \to \ZZ$ by 
$$
\varepsilon_0(a)=\sum_{\substack{a_i \equiv a \\ \mathrm{mod} \ d}} [K_i \colon K] \chi(D_i^{(\lambda), \circ}). 
$$

From the definition of the gamma factor, Remark \ref{gam-brie} and the expression for the residues in Theorem \ref{esn-vie}, we get:
$$
\Gamma(\nabla^{(d-\lambda)} \colon \Lcal^{(d-\lambda)^{-1}})=\prod_{i \in I} \Gamma \left(1-\{ \tfrac{a_i \lambda}{d}\} \right)^{[K_i \colon K] \chi(D_i^{(\lambda), \circ})} =\prod_{a=1}^{d-1} \Gamma\left(1-\frac{a}{d} \right)^{\varepsilon_0(a \slash \lambda)}
$$ and, putting everything together, 
$$
\prod_{j=0}^{2n} P_\lambda(\mathrm{det}_{\QQ(\xi)} H^j(Y, g))^{(-1)^j} \sim_{K^\times} \alpha_{\lambda} \cdot (2\pi i)^{M_\lambda} \cdot \prod_{a=1}^{d-1} \Gamma\left(1-\frac{a}{d}\right)^{\varepsilon_0(a /\lambda)}. 
$$ 

\begin{lemma}\label{thepair} We have $\alpha_\lambda^w \in K^\times$. 

\end{lemma}

\begin{proof} It suffices to show that the equality $(\Mca_\lambda, c_{Z \mathrm{mod} D^{(\lambda)}})^d=1$ holds in $\CC^\times \slash K^\times$. Since $c_{Z \mathrm{mod} D^{(\lambda)}}$ is a $\ZZ$-linear combination of closed points of $U^{(\lambda)}$, we are reduced to prove that $N_{K(x) \slash K}(\Mca_{\lambda | x})^d=1$ for any such $x$. We shall use the explicit description given in \eqref{explicit-pairing}. Consider the cyclic field extension $K(Y) \slash K(Z)$ corresponding to $\pi \colon Y \to Z$. By Kummer theory, there exists a generator $f_\lambda$ such that $f_\lambda^g=\xi^\lambda f_\lambda$. Since $x$ does not belong to the ramification locus, we may choose $f_\lambda$ invertible at $\pi^{-1}(x)$. Using part \ref{EVa} of Theorem \ref{esn-vie}, we see that $f_\lambda$ is a local basis of the invertible sheaf $\Lcal^{(\lambda)^{-1}}$ at $x$, hence we have $\varphi=f_\lambda^{-1}$, seen as a function on a small neighborhood of $x$ in the analytic manifold associated to $U^{(\lambda)}$. Therefore,
$$
N_{K(x) \slash K}(\Mca_{\lambda | x})^d=\prod_{\sigma \in \mathrm{Aut}_K(K(x))} f_\lambda(x^\sigma)^d=N_{K(x) \slash K}(f_\lambda^d(x))=1, 
$$ where the last equality holds because $f_\lambda^d$ is a rational function on $Z$. 
\end{proof}

\subsubsection*{Step 2} To conclude the proof, it remains to show:  

\begin{proposition}\label{chicar} The following equality holds: 
$$
\sum_{j=0}^{2n} (-1)^j \sum_{p+q=j} p h^{p, q}_\lambda(Y)-\frac{1}{2}\sum_{m \geq 0} (-1)^m (n+m) \chi(D^{(\lambda)}(m))=\sum_{i \in I} \{ \tfrac{a_i \lambda}{d} \} [K_i \colon K] \chi(D_i^{\lambda, \circ}). 
$$
\end{proposition}

We shall need several lemmas. 

\begin{lemma}\label{hdg} Let $X$ be a smooth projective variety and $D$ a simple normal crossings divisor on $X$, with everything defined over a subfield of $\CC$. Then: 
$$
\frac{1}{2} \sum_{m \geq 0} (-1)^m (n+m) \chi(D(m))=\sum_{p=0}^{n} (-1)^p p \chi(X, \Omega^p_X(\log D)). 
$$
\end{lemma}

\begin{proof} Recall that the complex $\Omega_X^\bullet(\log D)$ is equipped with a weight filtration $W$ and that the Poincar\'e residue induces an isomorphism 
$$
\Gr^W_m \Omega^\bullet_X(\log D) \stackrel{\sim}{\to} (a_m)_\ast \Omega^\bullet_{D(m)}[-m],
$$ where $a_m \colon D(m) \to X$ is the map deduced from the inclusions $D_J \hookrightarrow X$. From this, we compute 
\begin{align*}
\sum_{p=0}^n (-1)^p p \chi(X, \Omega^p_X(\log D))&=\sum_{m \geq 0} \sum_{p=0}^n (-1)^p p \chi(X, \Gr_m^W \Omega^p_X(\log D)) \\
&=\sum_{m \geq 0} (-1)^m \sum_{p=0}^{2(n-m)} (-1)^p (p+m) \chi(D(m), \Omega^p_{D(m)}) \\
&=\frac{1}{2} \sum_{m \geq 0} (-1)^m (n+m) \chi(D(m)). 
\end{align*} In what precedes, the last equality follows from 
$$
\frac{n-m}{2}\chi(D(m))=\sum_{p=0}^{n-m} (-1)^p p \chi(D(m), \Omega^p_{D(m)}), 
$$ which, regrouping terms, is an immediate consequence of Serre's duality. 
\end{proof}

\begin{corollary}\label{hdg-comp} The following identity holds: 
$$
\frac{1}{2} \sum_{m \geq 0} (-1)^m (n-m) \chi(D(m))=\sum_{p=0}^{n} (-1)^p p \chi(X, \Omega^p_X(\log D)(-D)). 
$$
\end{corollary}

\begin{proof} Since $\Omega^n_Z$ is isomorphic to $\Omega^n_Z(\log D)(-D)$, a direct application of Serre's duality shows that $h^q(Z, \Omega^p_Z(\log D)(-D))=h^{n-q}(Z, \Omega^{n-p}_Z(\log D)).$ Now the claim follows from Lemma \ref{hdg} by elementary manipulations. 
\end{proof}

\begin{lemma}\label{st-ch} Let $\Fca$ be a locally free $\Oc_Z$-module of rank $r$. Then: 
$$
\sum_{p=0}^r (-1)^p p \ch^p(\Lambda^p \Fca)=(-1)^r c_{r-1}(\Fca)+\text{higher order terms}. 
$$
\end{lemma}

\begin{proof}
See the proof of \cite[Lemma 5.2]{ST97}. 
\end{proof}

\begin{lemma}\label{1stch} The first Chern class of $\Lcal^{(\lambda)^{-1}}$ is $c_1(\Lcal^{(\lambda)^{-1}})=-\sum_{i \in I} \{\tfrac{a_i \lambda}{d}\} [D_i].$
\end{lemma}

\begin{proof} This is a straightforward consequence\footnote{Alternatively, it follows from the general expression of Chern classes of vector bundles with logarithmic connections in terms of its residues \cite[App. B]{EV86}.} of the definition of $\Lcal^{(\lambda)}$, together with the assumption that $\Lcal^d \simeq \Oc_Z(D)$. 
\end{proof}

\begin{lemma}\label{inter} Let $D$ be a simple normal crossings divisor on a smooth projective variety $X$, with everything defined over $K$. Then: 
$$
(-1)^{n-1}\deg_K c_{n-1}(\Omega^1_X(\log D)_{| D_i})=[K_i \colon K] \chi(D_i^\circ).  
$$
\end{lemma}

\begin{proof} Let us first recall that 
\begin{equation}\label{chern-log}
\chi(X \setminus D)=(-1)^n \deg_K c_n(\Omega^1_X(\log D)),  
\end{equation} an identity which follows from $\chi(X)=(-1)^n \deg_K c_n(\Omega^1_X)$ \cite[Ex. 8.1.12]{Ful84}, the exact sequences of logarithmic differentials \cite[\S 2.3]{EV92} and the additivity of the Euler characteristic. By the exact sequence of sheaves on $D_i$ 
$$
0 \to \Omega^1_{D_i}(\log(D-D_i)_{| D_i} ) \longrightarrow \Omega^1_X(\log D)_{| D_i} \longrightarrow \Oc_{D_i} \longrightarrow 0 
$$ and the Whitney sum formula for Chern classes, we get 
$$
c_{n-1}(\Omega^1_X(\log D)_{| D_i})=c_{n-1}(\Omega^1_{D_i}(\log (D-D_i)_{| D_i}). 
$$ Since $K_i$ is the field of constants of $D_i$, the claim follows from \eqref{chern-log} applied to the simple normal crossings divisor $(D-D_i)_{|D_i}=\bigcup_{i \neq j} (D_i \cap D_j)$ on $D_i$. 
\end{proof}

Armed with all these lemmas, we can now conclude the proof of Theorem \ref{maintheoD}. 

\begin{proof}[Proof of Proposition \ref{chicar}] Let $(\star)$ denote the right hand side of the equality to be proved. By part \ref{EVb} of Theorem \ref{esn-vie} and Lemma \ref{hdg}, it can be rewritten as 
$$
(\star)=\sum_{p=0}^n (-1)^p p [\chi(Z, \Lcal^{(\lambda)^{-1}} \otimes_{\Oc_Z} \Omega^p_Z(\log D^{(\lambda)})))-\chi(Z, \Omega^p_Z(\log D^{(\lambda)})) ]. 
$$

The Hirzebruch-Riemann-Roch theorem \cite[Cor. 15.2.1]{Ful84} and the multiplicativity of the Chern character yields the identity
\begin{align*}
\chi(Z, \Lcal^{(\lambda)^{-1}} \otimes_{\Oc_Z} \Omega^1_Z&(\log D^{(\lambda)}))-\chi(Z, \Omega^p_Z(\log D^{(\lambda)})) \\
&=\int_Z [\ch(\Lcal^{(\lambda)^{-1}})-1] \cdot \ch(\Omega^p_Z(\log D^{(\lambda)})) \cdot \Td(TZ).  
\end{align*} Therefore, 
\begin{align*}
(\star)&=\int_Z [\ch(\Lcal^{(\lambda)^{-1}})-1] \cdot [\sum_{p=0}^n (-1)^p p \ch(\Lambda^p \Omega^1_Z(\log D^{(\lambda)}))]\cdot \Td(TZ) \\
&=(-1)^n \int_Z c_1(\Lcal^{(\lambda)^{-1}}) c_{n-1}(\Omega^1_Z(\log D^{(\lambda)})) \\
&=\sum_{i \in I} \{\tfrac{a_i \lambda}{d} \} (-1)^{n-1} \deg_K c_{n-1}(\Omega^1_Z(\log D^{(\lambda)})_{| D_i^{(\lambda)}}) \\
&=\sum_{i \in I} \{\tfrac{a_i \lambda}{d} \} [K_i \colon K] \chi(D_i^{(\lambda), \circ}), 
\end{align*} where the last three equalities follow from lemmas \ref{st-ch}, \ref{1stch} and \ref{inter}.
\end{proof}

\subsection{Proof of Proposition \ref{freeaction}}\label{proof-free} \

\noindent 
Let us keep the notation from paragraph \ref{dethdg}, and denote by $\pi \colon X \to Z$ the quotient of $X$ by the action of $g$. Since $g$ has no fixed points, $Z$ is a smooth projective variety and the morphism $\pi$ is \'etale. Thus the direct image of the exterior differential decomposes into rank one \textit{smooth} connections, and the Saito-Terasoma theorem reduces to 
$$
\per(\Mca_\lambda)=\per(\unit_Z) \cdot (\Mca_\lambda, c_Z).  
$$ 

Let us compute each of these terms. On the one hand, by the identity \eqref{drB} and the multiplicativity of the Euler characteristic on \'etale covers, we have 
$$
\per(\unit_Z)=\sqrt{\disc(q_{dR}) \cdot \disc(q_{B})} \cdot (2\pi i)^{\frac{n \chi(X)}{2d}}, 
$$ where $q_{dR}$ and $q_B$ denote the cup-product on $H^n_{dR}(Z \slash K)$ and $H^n_B(Z)$ respectively. On the other hand, the canonical cycle is given by $c_Z=(-1)^n \mathrm{div}(\omega)$ for any non-zero rational section $\omega$ of $\Omega^1_Z$. Let $f_\lambda$ be a generator of $K(X)$ over $K(Z)$ such that $f_\lambda^g=\xi^\lambda f_\lambda$ and that $f_\lambda$ is invertible at the points lying over the support of $c_Z$. Then, by the proof of Lemma \ref{thepair},  
$$
(\Mca_\lambda, c_Z)=\prod_{z \in Z} \prod_{\sigma \in \mathrm{Aut}_K(K(z))} f_\lambda(z^\sigma)^{(-1)^n \ord_z(\omega)}.
$$ 

Therefore, for the obvious choice of $\alpha_\lambda$, 
$$
\prod_{j=0}^{2n} P_\lambda(H^j(X, g))^{(-1)^j} \sim_{K^\times} \alpha_\lambda (2\pi i)^{\frac{n\chi(X)}{2d}}. 
$$

\begin{remark}\label{hodge-classes} Let us assume that $H^j(X, g)$ vanishes except in degree $j=n$. In this case, Proposition \ref{chicar} shows that 
$$
(-1)^n \sum_{p+q=n} p h^{p, q}_\lambda(X)=\frac{n\chi(X)}{2d}. 
$$ Hence the function $\rho_n$ in \eqref{func-det} is constant and $\det_{\QQ(\xi)} H^n(X, g)$ consists of Hodge classes, as we claimed after the statement of Proposition \ref{freeaction}. \end{remark}

\subsection{Proof of Theorem \ref{maintheoB}} \

\noindent 
Let us keep the notation from paragraph \ref{dethdg}. Taking into account that $H^j(X, g)$ vanishes for $j \neq 1$ and that the residue of $\nabla^{(\lambda)}$ at $s$ is equal to $\frac{\ord_s(F_\lambda)}{d}$, the Saito-Terasoma theorem, the computation of $\per(\unit_{Z \setminus S})$ in Example \ref{sqrcurves} and the explicit description of the pairing in Example \ref{pairing-curves} yield: 
$$
P_\lambda(\mathrm{det}_{\QQ(\xi)} H^1(X, g))=\alpha_\lambda (2\pi i)^{g(Z)-1+\sum_{s \in S} [K(s) \colon K]} \prod_{s \in S} \Gamma\left(1-\frac{\ord_s(F_\lambda)}{d} \right)^{-[K(s) \colon K]}.
$$ 

To transform this into the desired expression, it suffices to show that 
$$
(2\pi i)^{\sum_{s \in S} [K(s) \colon K]} \prod_{s \in S} \Gamma\left(1-\frac{\ord_s(F_\lambda)}{d} \right)^{-[K(s) \colon K]} \sim_{K^\times} \prod_{s \in S} \Gamma\left(\frac{\ord_s(F_\lambda)}{d}\right)^{[K(s) \colon K]}. 
$$ Since $K$ contains $\QQ(\xi)$, this follows from the reflection formula 
\begin{equation}\label{reflec}
\Gamma(z)\Gamma(1-z)=\frac{\pi}{\sin(\pi z)}=2\pi i \frac{e^{i \pi z}}{e^{2\pi i z}-1}, 
\end{equation} together with the remark that $\sum_{s \in S} [K(s) \colon K] \frac{\ord_s(F_\lambda)}{d}$ is an integer since, using Lemma \ref{1stch}, it is equal to minus the degree of the Chern class $c_1(\Lcal^{(\lambda)^{-1}}).$

\subsection{Periods with compact support} \

\noindent 
For later use, we compute the periods with compact support of the object $\Mca_\lambda$ over the subvariety $U \subseteq U^{(\lambda)}$. 

Let $\pi \colon Y \to Z$ be the cyclic cover ramified at $D$ constructed in paragraph \ref{cyc-EV} and set $E=\pi^{-1}(D)$ and $V=Y \setminus E_{red}$. We denote by $\Sigma$ the singular locus of $Y$, by $\iota \colon Y \setminus \Sigma \to \Sigma$ the inclusion, and define the complex of coherent sheaves 
$$
\widetilde{\Omega}^\bullet_Y(\log E)=\iota_\ast \Omega^\bullet_{Y \setminus \Sigma}(\log(E \setminus \Sigma)).
$$ By \cite[\S1.17]{Ste77} and \cite[Lemma 1.2]{Ara14}, there is a canonical isomorphism
$$
\mathbb{H}^\bullet(Y, \widetilde{\Omega}^{\bullet}_Y(\log E)) \otimes_K \CC \stackrel{\sim}{\longrightarrow} H^j_B(V) \otimes_{\QQ} \CC,
$$
and the eigenspace decomposition of $\pi_\ast \widetilde{\Omega}^p_Y(\log E)$ is given by
$$
\pi_\ast \widetilde{\Omega}^p(\log E)=\bigoplus_{\lambda=0}^{d-1} \Lcal^{(\lambda)^{-1}} \otimes_{\Oc_Z} \Omega^p_Z(\log D). 
$$

Arguing as after Proposition \ref{esn-vie}, we get
$$
\per_c(\Mca_{\lambda | U})=\prod_{j=0}^{2n} P_\lambda(\mathrm{det}_{\QQ(\xi)} H^j_c(V, g))^{(-1)^j}, 
$$ where $H^j_c(V, g)$ denotes the largest sub-vector space $H \subseteq H^j_{B,c}(V)$ such that $g^\ast$ acts on $H \otimes_{\QQ} \CC$ through primitive eigenvalues, and the determinant is taken with respect to the natural structure of $\QQ(\xi)$-vector space. Finally, let us define\footnote{Note that this notation is coherent with the definition of cohomology with compact support in paragraph \ref{stconsup}, since the extension $\Lcal^{(\lambda)^{-1}}(D^{(\lambda)}-D)$ is small by Remark \ref{small-ext} and Theorem \ref{esn-vie}.} 
$$
h^{p, q}_{\lambda, c}(V)=h^q(Z, \Lcal^{(\lambda)^{-1}}(D^{(\lambda)}) \otimes_{\Oc_Z} \Omega^p_Z(\log D)(-D)). 
$$ 

\begin{theorem}\label{thm32} There exists a function $\varepsilon \colon \ZZ/d \to \QQ$ satisfying the conditions:
\begin{enumerate}
\item[(a)] The following equality holds for all $\lambda \in (\ZZ/d)^\times$:
\begin{equation*}
\sum_{a \in \ZZ/d} \varepsilon(a) \{\tfrac{a\lambda}{d}\}=\sum_{j=0}^{2n} (-1)^j \sum_{p+q=j} p h^{p, q}_{\lambda, c}(V).
\end{equation*}
\item[(b)] For each $\lambda \in (\ZZ/d)^\times$, there exists an algebraic number $\alpha_\lambda$ such that 
\begin{equation*}
\prod_{j=0}^{2n} P_\lambda(\mathrm{det}_{\QQ(\xi)} H^j_{c}(V, g))^{(-1)^j}\sim_{L^\times}  \alpha_\lambda \cdot \prod_{a=1}^{d-1} \Gamma\left(1-\frac{a}{d} \right)^{\varepsilon(a \slash \lambda)}, 
\end{equation*} where $L$ is the extension of $K$ generated by $i$ and $d^{\frac{1}{d-1}}$. Moreover, $\alpha_\lambda^w \in K^\times$. 
\end{enumerate}
\end{theorem}

\begin{proof} By Theorem \ref{st-compsup}, the following equality holds in $\CC^\times \slash K^\times:$
$$
\per_c(\Mca_{\lambda | U})=\per_c(\unit_U) \cdot (\Mca_{\lambda|U}, c_{Z \mathrm{mod} D}) \cdot \Gamma(\nabla^{(\lambda)} \colon \Lcal^{(\lambda)^{-1}})^{-1}. 
$$ 

The first factor was computed in corollaries \ref{prop-unit-cs} and \ref{hdg-comp}: $\per_c(\unit_U)=C (2\pi i)^M$, where $M=\sum_{p=0}^n (-1)^p p \chi(Z, \Omega^p_Z(\log D)(-D))$ and $C^2 \in K^\times$. Besides, the number $\alpha_\lambda=C (\Mca_{\lambda | U}, c_{Z \mathrm{mod} D})$ is algebraic and satisfies $\alpha_\lambda^w \in K^\times$ by Lemma \ref{thepair}. Let $I^{(\lambda)}$ denote the subset of $I$ consisting of those indices such that $a_i \lambda$ does not divide $d$ and define a function $\varepsilon_0 \colon \ZZ/d \to \ZZ$ by 
$$
\varepsilon_0(a)=\sum_{\substack{i \in I^{(\lambda)} \\ a_i \equiv a \ \mathrm{mod} \ d}} [K_i \colon K] \chi(D_i^\circ).  
$$ Combining the definition of the gamma factor with the reflection formula \eqref{reflec} and the remark right after it, as well as the multiplication formula \eqref{mult-gamma}, we get 
\begin{align*}
\per_c(\Mca_{\lambda | U})&\sim_{K^\times} \alpha_\lambda (2\pi i)^M \prod_{i \in I^{(\lambda)}} \Gamma \left(\{\tfrac{a_i \lambda}{d}\} \right)^{-[K_i \colon K] \chi(D_i^\circ)} \\
&\sim_{K^\times} \alpha_\lambda (2\pi i)^{M-\sum_{i \in I^{(\lambda)}} [K_i \colon K] \chi(D_i^\circ) } \prod_{a=1}^{d-1} \Gamma\left(1-\frac{a}{d} \right)^{\varepsilon_0(a \slash \lambda)} \\
&\sim_{L^\times} \alpha_\lambda \prod_{a=1}^{d-1} \Gamma\left(1-\frac{a}{d} \right)^{\varepsilon(a \slash \lambda)}, 
\end{align*} where $\varepsilon(a)=\varepsilon_0(a)+\frac{2}{d-1}M-\frac{2}{d-1} \sum_{i \in I^{(\lambda)}} [K_i \colon K] \chi(D_i^\circ)$. 

To conclude, it remains to show that 
\begin{align*}
\sum_{p=0}^n (-1)^p &p \chi(Z, \Lcal^{(\lambda)^{-1}}(D^{(\lambda)}) \otimes_{\Oc_Z}\Omega^p_Z(\log D)(-D)) \\ &-\sum_{p=0}^n (-1)^p p \chi(Z, \Omega^p_Z(\log D)(-D))=\sum_{i \in I^{(\lambda)}} (\{\tfrac{a_i \lambda}{d} \}-1)[K_i \colon K] \chi(D_i^\circ).  
\end{align*} \textit{Mutatis mutandi}, the proof is the same as the one of Proposition \ref{chicar}. 
\end{proof}

\section{Proof of Theorem \ref{maintheoA}}\label{sec-4}

In this final section, we show how to derive Theorem \ref{maintheoA} from the previous computations of periods of cyclic covers. As already mentioned, the strategy consists of performing a series of equivariant blow-ups to make $g$ act freely outside a simple normal crossings divisor. 

\subsection{The cohomology of equivariant blow-ups} \

\noindent
Let $X$ be a smooth projective variety over a subfield $K$ of $\CC$, together with the action of a finite group $G$, and let $Y \hookrightarrow X$ be a $G$-equivariant smooth closed subscheme of codimension $r$. Denote by $\tau \colon \widetilde{X} \to X$ the blow-up of $X$ along $Y$, by $e \colon E \hookrightarrow \widetilde{X}$ the immersion of the exceptional divisor and by $\psi \colon E \to Y$ the restriction of $\tau$ to $E$. Finally, let $\Oc_E(1)$ be the tautological vector bundle on $E$. By the universal property of blow-ups, the action of $G$ canonically lifts to $\widetilde{X}$ and the morphisms $\tau$ and $e$ are $G$-equivariant (see e.g. \cite[\S2.1.11]{IT14}).

\begin{lemma}\label{blowblow} The map 
\begin{align}\label{equiv-blow}
H^j(X) \oplus \bigoplus_{\ell=0}^{r-2} H^{j-2\ell-2}(Y)(-\ell-1) &\longrightarrow H^j(\widetilde{X}) \nonumber \\
(\eta, \kappa_0, \ldots, \kappa_{r-2}) &\longmapsto \tau^\ast \eta + e_\ast [\sum_{\ell=0}^{r-2} \psi^\ast(\kappa_\ell) \cup c_1(\Oc_E(1))^\ell ]
\end{align} is a $G$-equivariant isomorphism in the category $M_{K, F}(K)$. Moreover, it is compatible with the Hodge structures on the Betti components. 
\end{lemma}

\begin{proof} Leaving aside the action of $G$, the statement is proved in \cite[\S 9]{Man68}. To show that \eqref{equiv-blow} is $G$-equivariant, it suffices to observe that $G$ acts trivially on the Chern class $c_1(\Oc_E(1))$, which follows from the fact that the vector bundle $\Oc_E(1)$ is $G$-linearized. Indeed, $G$ acts on the normal bundle $N_Y X$ because the subscheme $Y$ is $G$-stable, hence on the geometric vector bundle associated to $N_Y X$ and on its blow-up along the zero section. But this blow-up is nothing but the geometric vector bundle corresponding to $\Oc_E(1)$, which yields the $G$-linearization. 
\end{proof}

For the rest of the section, we assume that $G$ is cyclic of order $d \geq 3$ and we choose a generator $g$. From Lemma \ref{blowblow} we immediately deduce: 

\begin{corollary}\label{cor-act-bl} If the restriction of $g$ to $Y$ has order strictly less than $d$, then the Hodge structures $H^j(X, g)$ and $H^j(\widetilde{X}, g)$ are isomorphic. 
\end{corollary}

We shall perform a sequence of blow-ups with the goal of making $G$ act freely outside a $G$-strict normal crossings divisor $D$. Recall that this means that, for each irreducible component $D_i$ of $D$, either $gD_i=D_i$ or $gD_i \cap D_i=\emptyset$. Given a subgroup $H$ of $G$, let $X^H$ denote the closed subscheme of fixed points of $H$. We shall make repeated use of the fact that $X^H$ is smooth, see e.g. \cite[Prop. 3.1]{Edi92}. 

Let $d=p_1^{e_1} \cdots p_t^{e_t}$ be the prime factorization of $d$. Put $N=\sum_{i=1}^t e_t$ and let $S$ denote the set of divisors of $d$. We can write $S$ as a disjoint union of subsets 
$$
S_\ell=\{ p_1^{f_1} \cdots p_t^{f_t} \ | \ 0 \leq f_i \leq e_i, \ f_1+\cdots+f_t=\ell \}
$$ for $\ell=0, \ldots, N$. For example, $S_0=\{1\}$ and $S_1=\{p_1, \ldots, p_t\}$. Given an integer $s \in S$, let $H_s$ denote the subgroup of $G$ generated by $g^s$. 

We construct a tower of blow-ups
$$
\widetilde{X}=X_{N-1} \xrightarrow{\tau_{N-1}} \cdots \xrightarrow{\tau_1} X_0 \xrightarrow{\tau_0} X
$$ as follows: $\tau_0 \colon X_0 \to X$ is the blow-up of $X$ along $X^G$ and, for each $m \geq 1$, we define recursively $\tau_m \colon X_m \to X_{m-1}$ as the blow-up of $X_{m-1}$ along the strict transform of $\bigcup_{s \in S_m} X^{H_s}$ by $\tau_{m-1} \circ \cdots \circ \tau_0$. We denote by $\tau \colon \widetilde{X} \to X$ the composition of $\tau_{N-1} \circ \cdots \circ \tau_0$ and by $E$ the inverse image of $\bigcup_{s \in S} X^{H_n}$ by $\tau$. 

Since all the centers of the blow-ups are smooth and $G$-equivariant, $\widetilde{X}$ is a smooth projective variety over $K$, acted upon by an automorphism of order $d$ which will be still denoted by $g$. 

\begin{lemma}\label{blows}
\begin{enumerate}
\item[]
\item[\namedlabel{BLa}{(a)}]  $E$ is a $G$-strict normal crossings divisor and $G$ acts freely on $\widetilde{X}\setminus E$. 
\item[\namedlabel{BLb}{(b)}] $G$ acts with order strictly less than $d$ on the cohomology\footnote{Note, however, that $G$ does not need to act with order strictly less than $d$ on $E$ itself; the action may not be trivial even when one only blows up fixed points.} of $E$. 
\item[\namedlabel{BLc}{(c)}] The Hodge structures $H^j(X, g)$ and $H^j(\widetilde{X}, g)$ are isomorphic. 
\end{enumerate}
\end{lemma}

\begin{proof} Assertion \ref{BLa} is proved in \cite[Lemma 4.2.3]{IT14}. Property \ref{BLb} follows from the fact that $\widetilde{X}$ is obtained by successively blowing up $G$-equivariant smooth closed subschemes on which $G$ acts with order strictly less than $d$. Indeed, at each step, the exceptional divisor is a projective bundle over the center $Y$ of the blow-up, so its cohomology ring is an algebra over $H^\ast(Y)$ generated by the first Chern class of the tautological line bundle. By the proof of Lemma \ref{blowblow} above, $G$ acts trivially on this Chern class, so the order of the action is the same as the order of the action on $H^\ast(Y)$, which is strictly smaller than $d$ by construction. Finally, \ref{BLc} is a straighforward consequence of Corollary \ref{cor-act-bl}. 
\end{proof}

\subsection{End of the proof} \

\noindent
Let us keep the notation of Theorem \ref{maintheoA}, and write $G$ for the group generated by the automorphism $g$. By Lemma \ref{blows}, we may assume that there exists a simple normal crossings divisor $E$ on $X$ such that $G$ acts freely on the complement $V=X \setminus E$ and with order strictly less than $d$ on the cohomology of $E$. Since the diagram 
\begin{equation}\label{excision}
\xymatrix{H^{j-1}_{dR}(E) \otimes \CC \ar[r] \ar[d] & H^j_{dR, c}(V) \otimes \CC \ar[r] \ar[d] & H^j_{dR}(X) \otimes \CC \ar[r] \ar[d] & H^j_{dR}(E) \otimes \CC \ar[d] \\
H^{j-1}_B(E) \otimes \CC \ar[r] & H^j_{B, c}(V) \otimes \CC \ar[r] & H^j_B(X) \otimes \CC \ar[r] & H_B^j(E) \otimes \CC 
}, 
\end{equation} commutes and all the maps are $G$-equivariant, we deduce that
$$
P_\lambda(\mathrm{det}_{\QQ(\xi)} H^j(X, g))=P_\lambda(\mathrm{det}_{\QQ(\xi)} H^j_c(V, g)).
$$

Let $U$ denote the quotient $V \slash G$. Since $G$ acts freely, $U$ is a smooth quasi-projective variety over $K$ and $V$ is an unramified cyclic cover of $U$. Let $Z$ be a smooth projective variety over $K$ containing $U$ as the complement of a normal crossings divisor with irreducible components $(D_i)_{i \in I}$, and consider the normalization $Y$ of $Z$ in the fraction field of $V$. Then the action of $G$ extends to $Y$, and the normalization map $\pi \colon Y \to Z$ exhibits $Z$ as the quotient of $Y$ by $G$.

\begin{lemma} $\pi \colon Y \to Z$ is isomorphic to a cyclic cover constructed out of an invertible sheaf $\Lcal$ on $Z$ such that $\Lcal^{d} \simeq \mathcal{O}_Z(D)$, for some effective divisor $D$ supported on $\bigcup_{i \in I} D_i$.
\end{lemma}

\begin{proof} Since $Z$ is smooth and $Y$ is an abelian normal cover, the morphism $\pi$ is flat and there is a direct sum decomposition
$$
\pi_\ast \mathcal{O}_Y=\bigoplus_{\chi \in \widehat{G}} L_\chi^{-1}, 
$$ where $L_\chi^{-1}$ denotes the invertible subsheaf of $\pi_\ast \mathcal{O}_Y$ on which $G$ acts through the character $\chi$ \cite[\S 3]{Ber84}. We claim that $\mathcal{L}=L_{\chi}$, for any primitive $\chi$, has the desired property. Indeed, if $f_\chi$ is a generator of the cyclic extension $K(Y) \slash K(Z)$ on which $G$ acts through $\chi$, then $g_\chi=f_\chi^d$ is a rational function on $Z$ whose divisor can be uniquely written as $\mathrm{div}(g_\chi)=-d B_\chi+D_\chi$, where $D_\chi$ is effective, supported on $\bigcup_{i \in I} D_i$ and has no multiplicity $\geq d$. If we pick another generator, $D_\chi$ remains unchanged, while $B_\chi$ is replaced by a linearly equivalent divisor. Therefore, $D_\chi$ and the invertible sheaf $\mathcal{O}_Z(B_\chi)$ are uniquely determined by $\chi$. A local computation, as in \textit{loc.cit.}, shows that $L_\chi=\mathcal{O}_Z(B_\chi)$, hence $\mathcal{L}^d=\mathcal{O}_Z(D_\chi)$.  
\end{proof}

Once that we have proved that $\pi \colon Y \to Z$ is a cyclic cover in the sense of paragraph \ref{cyc-EV}, the statement follows from Theorem \ref{thm32}, together with the observation that $h^{p, q}_\lambda(X)=h^{p, q}_{\lambda, c}(V)$ for all $\lambda \in (\ZZ/d)^\times$, since the lower horizontal maps in \eqref{excision} are morphisms of mixed Hodge structures. This concludes the proof. 

\medskip

\begin{small}
\begin{ack} The results of this paper were obtained during my PhD thesis under the supervision of C. Soul\'e and J. Wildeshaus, whom I would like to thank for their help and encouragement. In my long list of debts, S, Bloch occupies a special place, since his was the intuition that periods of CM motives should be related to epsilon factors of connections; this was the subject of a letter to H. Esnault \cite{Blo05}. I am grateful to her for generously sharing her insights in several conversations. Many thanks also to P. Deligne for two wonderful \textit{apr\`{e}s-midi} at the IH\'ES. The comments of J-B. Bost and T. Saito to preliminary versions of this paper helped me improve the expositions and correct some inaccuracies. Finally, I would like to thank G. Ancona, O. Benoist, J. I. Burgos Gil, P. Colmez, B, Edixhoven, D. Eriksson, G. Freixas i Montplet, B. H. Gross, M. Maculan, D. R\"ossler and C. Sabbah for useful discussions. 
\end{ack}
\end{small}

\end{document}